\documentclass{article}
\usepackage[utf8]{inputenc} 
\usepackage[T1]{fontenc}    
\usepackage{tikz}
\usepackage{algorithm}
\usepackage{algpseudocode}
\usepackage{float}
\usepackage{amssymb}
\usepackage{amsmath}
\usepackage{hyperref}
\usepackage{mathtools}
\usepackage{graphicx}
\usepackage{subcaption}
\usepackage{siunitx}
\usepackage{booktabs}
\usepackage{tikz-cd}
\usepackage{stix}
\usepackage{bm}
\usepackage{amsbsy}
\usepackage{amsthm}
\usepackage{arxiv}

\newtheorem{theorem}{Theorem}
\newtheorem{corollary}{Corollary}
\newtheorem{remark}{Remark}

\newcommand{\velocity}{\bm{v}}

\newcommand{\stdbase}{{\bm{e}}}

\newcommand\norm[1]{\left\lVert#1\right\rVert}
\DeclareMathOperator*{\argmax}{arg\,max}
\DeclareMathOperator*{\argmin}{arg\,min}

\newcommand{\inflowBoundary}{\Gamma_{-}}
\newcommand{\characteristicBoundary}{\Gamma_{0}}
\newcommand{\outflowBoundary}{\Gamma_{+}}

\newcommand{\noise}{\varepsilon}

\newcommand{\sourceIntensity}{\pmb{\lambda}}
\newcommand{\sourceIntensityi}[1]{\ensuremath{\lambda_{#1}}}

\newcommand{\misfit}{\bm{y}}
\newcommand{\misfiti[1]}{{y_{#1}}}

\newcommand{\measurement}{\bm{d}}
\newcommand{\measurementPlain}{\bm{d}}

\newcommand{\numberOfSensors}{N_d}

\newcommand{\concentration}{u}

\newcommand{\normal}{\bm{n}}
\newcommand{\wind}{{\bm{v}}}

\newcommand{\obsO}{\mathcal{B}} 
\newcommand{\pto}{\mathcal{F}} 
\newcommand{\pts}{\mathcal{K}} 
\newcommand{\som}{\mathcal{M}^+} 
\newcommand{\somN}{\mathcal{M}^+_N}
\newcommand{\sop}{D} 
\newcommand{\bochner}{\mathcal{V}} 
\newcommand{\mta}{\mathcal{Q}} 
\newcommand{\adjoint}{q} 
\newcommand{\Nobservations}{N_{\text{obs}}} 

\newcommand{\ndof}{{n_{\text{dof}}}}

\newcommand{\measure}{\mu}

\newcommand{\measurementSpace}{\mathbb{R}^{\Nobservations}}

\newcommand{\R}{\mathbb{R}}
\newcommand{\N}{\mathbb{N}}
\newcommand{\parameterContinuous}{m_{\text{C}}}
\newcommand{\parameterInitial}{m_{\text{I}}}

\newcommand{\measureInitial}{\measure_{\text{I}}}
\newcommand{\measureContinuous}{\measure_{\text{C}}}
\newcommand{\bmeasureInitial}{\bar{\measure}_{\text{I}}}
\newcommand{\bmeasureContinuous}{\bar{\measure}_{\text{C}}}

\newcommand{\ansatzSources}{\mathcal{S}} 
\newcommand{\ansatzSourcesInitial}{\mathcal{S}_{\text{I}}} 
\newcommand{\ansatzSourcesContinuous}{\mathcal{S}_{\text{C}}}

\newcommand{\observationforeqx}{x_i^\text{obs}}
\newcommand{\observationforeqt}{t^{\text{obs}}_i}
\newcommand{\observationforequation}{(\observationforeqt,\observationforeqx)}
\newcommand{\x}{x}

\newcommand{\xb}{\bm{x}}
\newcommand{\xbr}{{\bm{x}}}

\newcommand{\up}{{\vec{u}^{n+1}}}
\newcommand{\un}{{\vec{u}^{n}}}
\newcommand{\pp}{{\vec{\adjoint}^{n+1}}}
\newcommand{\pn}{{\vec{\adjoint}^{n}}}
\newcommand{\ansatzSpace}{\mathcal{V}_h}
\newcommand{\FEMi}{I}
\newcommand{\ellipI}{\eta}
\newcommand{\ellipLaplace}{\gamma}
\newcommand{\Rplus}{\mathbb{R}_{>0}}

\newcommand{\fracSolidus}[2]{#1/#2\,}

\renewcommand{\vec}{\bm}

\newcommand{\tcr}[1]{\textcolor{black}{#1}}

\title{Sparse Source Identification in Transient Advection-Diffusion Problems with a Primal-Dual-Active-Point Strategy}

\date{} 					

\author{ \href{https://orcid.org/0009-0009-6325-3578}{\includegraphics[scale=0.06]{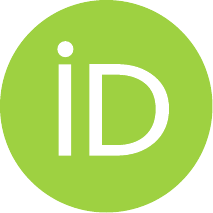}\hspace{1mm}Marco Mattuschka\textsuperscript{$1,$}}\thanks{Corresponding author e-mail: \texttt{marco.mattuschka@dlr.de}} ,\quad \href{https://orcid.org/0000-0003-3647-0728}{\includegraphics[scale=0.06]{orcid.pdf}\hspace{1mm}Daniel Walter\textsuperscript{$2$}},\quad
\href{https://orcid.org/0000-0002-2814-0027}{\includegraphics[scale=0.06]{orcid.pdf}\hspace{1mm}Max von Danwitz\textsuperscript{$1$}}, \quad
\href{https://orcid.org/0000-0002-8820-466X}{\includegraphics[scale=0.06]{orcid.pdf}\hspace{1mm}Alexander Popp\textsuperscript{$1,3$}}\\
\textsuperscript{$1$}German Aerospace Center (DLR), Institute for the Protection of Terrestrial Infrastructures, 53757 Sankt Augustin, Germany\\
\textsuperscript{$2$}Humboldt-Universität zu Berlin, Institut für Mathematik, 10117 Berlin, Germany\\
\textsuperscript{$3$}University of the Bundeswehr Munich, Institute for Mathematics and Computer-Based Simulation (IMCS),\\ 85577 Neubiberg, Germany
}

\renewcommand{\headeright}{}
\renewcommand{\undertitle}{}
\renewcommand{\shorttitle}{Sparse Source Identification in Advection-Diffusion Problems}

\begin{document}

\maketitle

\begin{abstract}
This work presents a mathematical model to enable rapid prediction of airborne contaminant transport based on scarce sensor measurements. The method is designed for applications in critical infrastructure protection (CIP), such as evacuation planning following contaminant release. In such scenarios, timely and reliable decision-making is essential, despite limited observation data. To identify contaminant sources, we formulate an inverse problem governed by an advection–diffusion equation. Given the problem's underdetermined nature, we further employ a variational regularization ansatz and model the unknown contaminant sources as distribution over the spatial domain. To efficiently solve the arising inverse problem, we employ a problem-specific variant of the Primal-Dual-Active-Point (PDAP) algorithm which efficiently approximates sparse minimizers of the inverse problem by alternating between greedy location updates and source intensity optimization. The approach is demonstrated on two- and three-dimensional test cases involving both instantaneous and continuous contaminant sources and outperforms state-of-the-art techniques with $L^2$-regularization. Its effectiveness is further illustrated in complex domains with real-world building geometries imported from OpenStreetMap.
\end{abstract}


\section{Introduction}
\begin{figure}
\begin{subfigure}{0.48\textwidth}
\includegraphics[width=\linewidth,trim={0 0 5cm 0},clip]{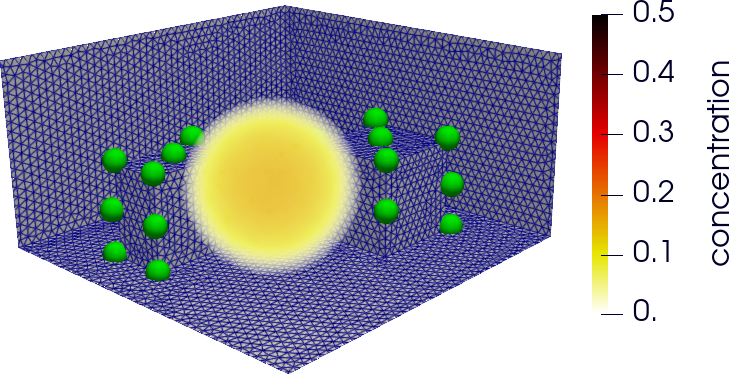}
\caption{}
\end{subfigure}
\begin{subfigure}{0.48\textwidth}
\includegraphics[width=\linewidth,trim={0 0 5cm 0},clip]{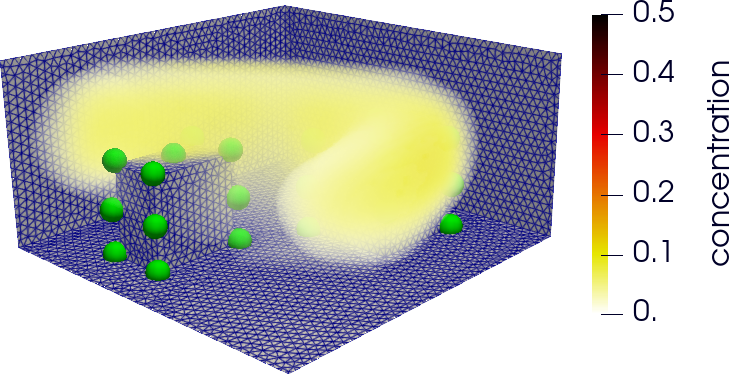}
\caption{}
\end{subfigure}
\begin{subfigure}{0.48\textwidth}
\includegraphics[width=\linewidth,trim={0 0 5cm 0},clip]{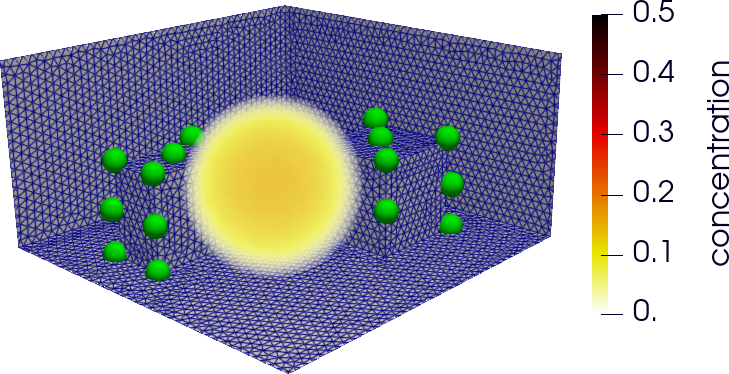}
\caption{}
\end{subfigure}
\begin{subfigure}{0.48\textwidth}
\includegraphics[width=\linewidth,trim={0 0 5cm 0},clip]{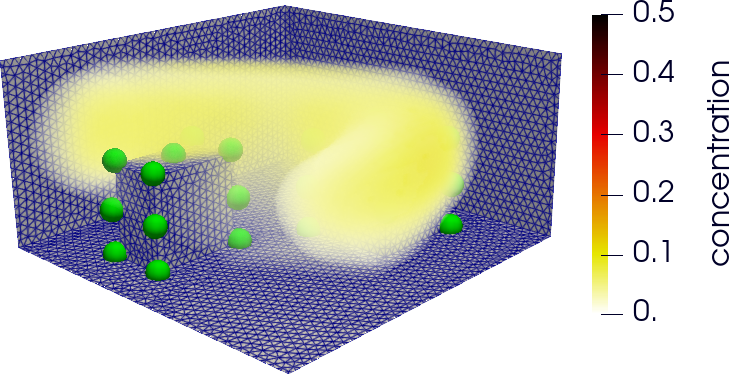}
\caption{}
\end{subfigure}
\centering
\begin{subfigure}{0.48\textwidth}
\centering
\includegraphics[width=0.5\linewidth]{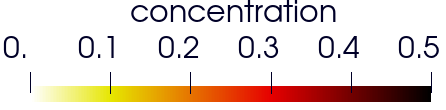}
\end{subfigure}
\caption{Three-dimensional gas source identification. Manufactured initial condition (a) and simulation at $t=\SI{3.0}{\second}$ (b), reproduced by algorithmic reconstruction of initial contaminant source (c) and simulation at $t=\SI{3.0}{\second}$ (d) based on sparse measurements. Positions of noisy point-evaluations (modeled sensors)
    are marked by green spheres.}
\label{fig:3D_Sim_Single}
\end{figure}
Airborne contaminant transport of hazardous materials poses a significant threat to communities and critical infrastructures. Contaminant release may occur accidentally in industrial leaks or spills, or may be caused intentionally in an act of sabotage or terrorism~\cite{Boris.2002, Patnaik.2012, Danwitz.2024}. An example visualization of the dispersion of a contaminant is provided in~\autoref{fig:3D_Sim_Single}. In order to stem these threats and predict further spreading of the contaminant, identifying the source of contamination is of utmost importance. In many scenarios, chemical hazards are caused by toxic industrial chemicals (TICs) including colorless gases invisible to the human eye and RGB cameras. Hence, inference about the unknown source is only possible indirectly through scarce measurements of the concentration field, e.g. point observations in the computational domain (green spheres in~\autoref{fig:3D_Sim_Single}). Describing the contaminant dispersion by an advection-diffusion equation and assuming models for the measurement process as well as the shape of the contaminant source, leads to an ill-posed inverse problem. In this paper, we revisit a recently popularized variational regularization ansatz for its stable solution and showcase its viability in scenarios with realistic computational domains.
Due to its practical relevance for decision-makers in emergency response, numerous approaches to source identification have been explored and documented in the literature. The specific computational method for source identification or source term estimation (STE) follows the underlying dispersion model. The complexity of contaminant dispersion models ranges from simple box models and Gaussian plume models to three-dimensional geometry-aware high fidelity models. An overview of commonly employed dispersion models is given in~\cite{Holmes.2006}. For an overview of the model-associated source term estimation methods, the interested reader is referred to~\cite{Hutchinson.2017} and references therein. The more specific task of active source identification additionally involves the routing of a mobile sensor system~\cite{Khodayimehr.2019}. Also, the very recent literature documents several practical approaches contributing to gas source localization~\cite{Hinsen.2024, Ruiz.2024, Wiedemann.2024}.

From the set of available dispersion models, we choose to model the airborne contaminant transport with a linear advection-diffusion equation in which the unknown sources appear as initial condition or right-hand-side term. In general, the passive transport of a substance (pollutant, chemical species, here: contaminant) in an incompressible or compressible fluid can be mathematically modeled by the advection-diffusion problem~\cite{Arystanbekova.2004, Ulfah.2018}. The wind vector field that drives the transport in the advection-diffusion problem might stem from Computational Fluid Dynamics (CFD) simulations~\cite{Maronga.2020, Danwitz.2019} or can be inferred from measurements itself~\cite{Towers.2016, Chen.2025}. In the present work, we consider the wind vector field as a known (fixed) input to the advection-diffusion equation. Throughout the paper, we refer to the process of simulating the latter for a given source and taking measurements of the obtained concentration field as the \emph{forward problem}. Given noisy measurements, the identification of the unknown source then corresponds to the associated ill-posed \emph{inverse problem}. 

For the particular application under consideration, various approaches for the stable solution of the inverse problem have been considered in the literature. For example, a Bayesian ansatz for finite-dimensional parametrized sources is presented in which the resulting posterior distribution is sampled via Markov chain Monte Carlo (MCMC) methods \cite{Albani.2021}. This approach heavily relies on the assumed measurement model which consists of a finite number of point measurements at a fixed time\tcr{,} allowing to compute all required PDE solves in an offline phase. In contrast, the present work considers time series data collected at a finite number of spatially fixed or moving sensor locations, thus preventing similar arguments. Instead, we formulate the source identification problem as a minimization task, incorporate the advection-diffusion equation as a constraint, and compensate the scarcity of measurement data by employing regularization techniques \cite{Gorelick.1983, Tsai.2014,  Villa.2021}. In this regard, Tikhonov regularization with weighted $L^2$-type penalty terms leads to quadratic minimization problems which are amenable to Newton methods, \cite{Wiedemann.2024, Petra.2011}. This approach can be further developed towards applications with real-time requirements by replacing high-fidelity models with reduced ones that still capture the essential characteristics, yet are much faster to evaluate~\cite{Key.2023, Dinkel.2024}. From a Bayesian perspective, problems of the described form appear in the computation of the maximum a posteriori estimator wherein specific forms of data misfit and penalty term are related to the noise model as well as to the prior distribution of the unknown source, respectively,~\cite{Stuart.2010,Nitzler.2022}. 

However, in realistic scenarios, it is reasonable to assume that contaminant sources are very localized, i.e. around an unknown finite number of locations, and that measurements are scarcely distributed in space, e.g. given by time series concentration data of few fixed sensors. Due to the smoothing behavior of the advection-diffusion process, weighted $L^2$-type regularization is unable to capture this structural assumption accurately. As a remedy, in the present paper, we assume that the sought-for contamination source is given by a superposition of finitely many simple atoms \tcr{of a given shape, e.g. point sources or Gaussians,} each characterized by a location parameter as well as a positive amplitude. Since their number is in general also unknown, this modeling assumption leads to a highly nonlinear parameter identification problem.   

\tcr{In order to alleviate these difficulties, we relax it from sums of atoms to parameters represented as integrals of the shape function w.r.t an unknown, positive Radon measure. The identification of the latter is tackled via variational regularization techniques in which the prior assumption on the sparse representation of the unknown source is incorporated by a Radon norm penalty. This approach replaces the combinatorial complexity and strong nonlinearity of the parametrized ansatz with a different set of challenges, namely, nonsmoothness and an infinite-dimensional space of unknowns, while offering the added benefit that the relaxed problem is convex. Moreover, it turns out that this relaxation is exact since optimal reconstructions consisting of finitely many Dirac-Delta functionals exist, see \autoref{thm:optimality}. Variational regularization approaches involving the Radon norm have received increased attention over the last years, see e.g. \cite{Bredies.2013,Scherzer.2009,Huynh.2024} for a nonexhaustive list of references. For the particular application to initial value identification in parabolic problems we refer to \cite{Casas.2019, Leykekhman.2020, Biccari.2023, Monge.2020}. In particular, several works have addressed the super-resolution properties of this ansatz, \cite{Fernandez.2014,Duval.2015}, i.e. the capability to recover the correct number of sources together with quantitative convergence guarantees for the associated points and amplitudes despite the presence of measurement noise. In this work, we adopt the framework outlined above for scenarios in contamination transport where measurements are scarcely distributed in space, e.g. given by time series concentration data of few fixed sensors. Our practical realization relies on the Primal-Dual-Active-Point (PDAP) strategy from \cite{Pieper.2021} which efficiently approximates sparse minimizers by alternating between greedy location updates and intensity optimization. Crucially, in comparison to e.g. \cite{Denoyelle.2020, Hnatiuk.2025,Chizat.2022}, this approach does not rely on nonconvex point-moving steps and is thus compatible with standard finite element discretizations. 
While, as outlined above, the concept of Radon norm regularization has been studied before, our particular focus lies on its application to practically relevant two and three-dimensional benchmark settings with real building imprints, bridging the gap between mathematical theory and the demand of real-world scenarios. Comparisons to state-of-the art methods with weighted $L^2$-type regularization, reveal excellent recovery results across a variety of examples while the proposed approach also exhibits a favorable computational scaling behavior w.r.t. to the availability of additional information, e.g., by adding new  sensors or higher sampling frequency.}

The remainder of the present paper is organized as follows. Mathematical modeling of airborne contaminant transport, along with the formulation of a sparsity promoting regularization approach for inverse problems, is presented in \autoref{sec:background}. For the sake of completeness, this \tcr{section} also briefly reviews a classical method with $L^2$-regularization. In~\autoref{sec:pdap}, the PDAP algorithm for the efficient solution the source identification problem is presented. Next, \autoref{sec:discretization} describes the finite element discretization of the forward and adjoint problem, followed by two- and three-dimensional numerical examples presented in~\autoref{sec:numerical}. Finally, \autoref{sec:conclusion} offers conclusions and an outlook to future developments.

\section{Background}\label{sec:background}
\subsection{Forward problem, parameter-to-observable and misfit-to-adjoint map}

\begin{figure}
\begin{subfigure}{0.48\textwidth}
\includegraphics[width=.9\linewidth]{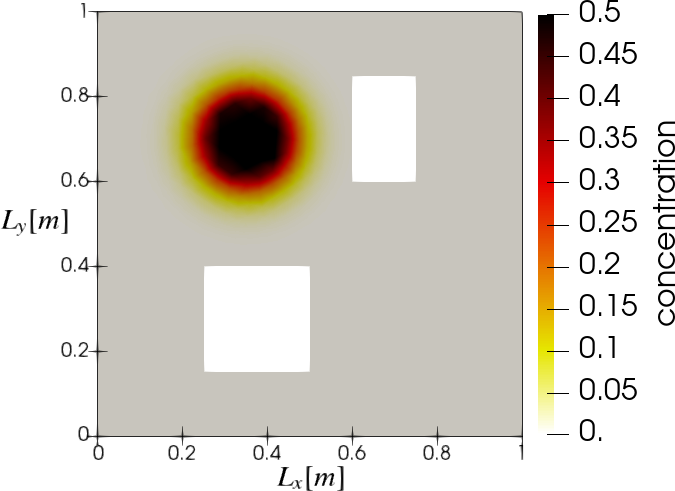}
\end{subfigure}
\begin{subfigure}{0.48\textwidth}
\includegraphics[width=.9\linewidth]{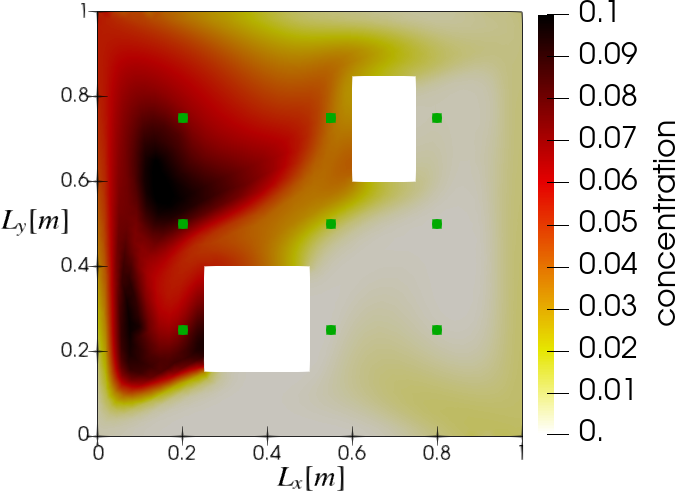}
\end{subfigure}
\caption{Illustration of commonly used benchmark problem~\cite{Villa.2021}. Ground truth initial condition $\parameterInitial$ (left), and ground truth concentration field at $t=\SI{5}{\s}$ (solution of \eqref{eq:forward_equation}), as well as measurement positions marked by nine green squares (right)}
\label{fig:example_tb_single}
\end{figure}
A mathematical description of the transport of a substance (contaminant) concentration $u$ in a bounded open domain $\Omega \subseteq \mathbb{R}^n \text{ for } n\in\{2,3\} $ is given by the following equation:

\begin{equation}\label{eq:forward_equation}
\begin{aligned}
u_t-\kappa\Delta u + \velocity\cdot\nabla u &= \parameterContinuous &\qquad&\text{in}\ (0,T)\times\Omega,\\
  \kappa\nabla u \cdot \normal &= 0 &&\text{in}\ (0,T)\times (\outflowBoundary \cup \characteristicBoundary),\\
  u&= 0 &&\text{in}\ (0,T)\times \inflowBoundary,\\
  u(0,\cdot) &= \parameterInitial &&\text{in}\ \Omega.
\end{aligned} \tag{$\mathcal{P_{\pts}}$}
\end{equation}
The problem~\eqref{eq:forward_equation} depends on two unknown functions $m=(\parameterInitial,\parameterContinuous)$ on $\Omega$, an initial condition $\parameterInitial$ and a source term $\parameterContinuous$ which we assume to be independent of time. The abstract space of admissible parameters will be denoted by $D$. A commonly used benchmark problem for initial value identification in advective-diffusive transport is illustrated in \autoref{fig:example_tb_single} (see also~\cite{Villa.2021,Petra.2011}). Note that the transport depends on a wind vector field, which is assumed to be smooth, bounded and divergence free, i.e., $\nabla \cdot \velocity = 0$. According to the relation of the wind vector and the outward-pointing boundary normal $\normal$, the domain boundary is subdivided into the outflow boundary $\outflowBoundary \subset \partial\Omega, \wind \cdot \normal > 0$, the inner boundary $\characteristicBoundary \subset \partial\Omega, \wind \cdot \normal = 0$ and the inflow boundary $\inflowBoundary \subset \partial\Omega, \wind \cdot \normal < 0$, comparable to \cite{Elman.2020}. \tcr{Finally, we emphasize that under suitable regularity assumptions on the problem data, i.e. the wind vector field $\velocity$ as well as the parameter space, $\eqref{eq:forward_equation}$ admits a unique solution $u$, in an appropriate weak sense, for given $m$, by standard parabolic regularity theory (see, e.g.,~\cite{Friedman.2008}). We note that significantly less regularity is required for the parameter $\parameterInitial$ compared to $\parameterContinuous$. For this reason, we may also admit Radon measures as initial conditions and, in particular, finite combinations of Dirac-Delta functionals.} 

\begin{figure}
\begin{subfigure}{0.48\textwidth}
\includegraphics[width=0.96\linewidth]{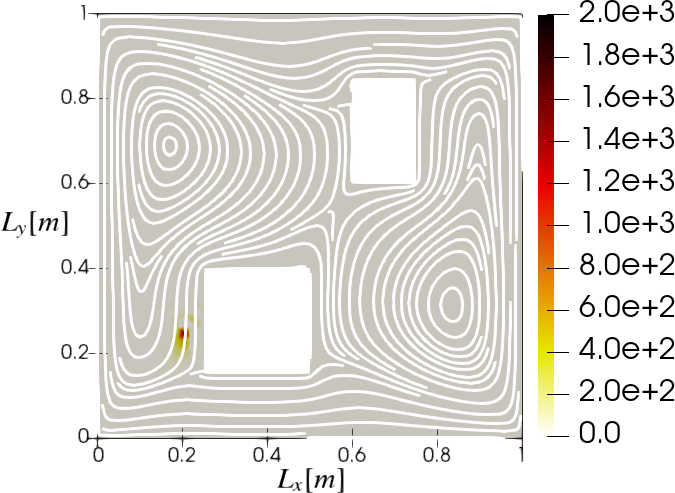}
\end{subfigure}
\begin{subfigure}{0.48\textwidth}
\includegraphics[width=0.90\linewidth]{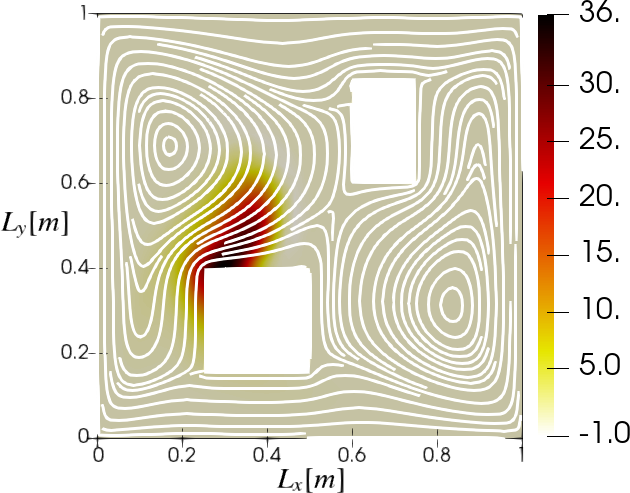}
\end{subfigure}
\caption{Adjoint transport problem. Example right-hand side of \eqref{eq:adjoint_equation}, i.e., a sensor misfit of ``$1$'' applied at $\observationforequation = (\SI{5}{\s}, [\SI{ 0.2}{\m} ,\SI{0.25}{\m}])$ (left), solution field of \eqref{eq:adjoint_equation} evaluated at $t=0$.}
\label{fig:adjoint}
\end{figure}

Next, we describe the observation setup. For this purpose, let $T_0>0$ be a fixed time point and consider a finite number of space-time points $\observationforequation$, $i=1,\dots,\Nobservations$, with $\observationforeqt \in [T_0,T)$ and $\observationforeqx \in \Omega$. The discrete space-time observation points are covered by open balls of radius \( r_0 > 0 \), defining the set
$\bar{\Omega}_o \coloneqq \bigcup_{i=1}^{\Nobservations} B_{r_0}(\observationforeqx) \subset \Omega.$
Choosing \( r_0 \) sufficiently small ensures that the solution to \eqref{eq:forward_equation} admits a continuous representative on \([T_0, T] \times \bar{\Omega}_o\). This gives rise to a linear and bounded \textit{space-time observation operator}
\begin{equation*}
\obsO\colon C([T_0, T] \times \bar{\Omega}_o) \to \mathbb{R}^{\Nobservations}, \quad   u \mapsto \sum_{i=1}^{\Nobservations}\delta_{\observationforequation}(u)\,\stdbase_i = \left(u\observationforequation\right)_{i=1}^{\Nobservations} 
\end{equation*}
where $\delta_{\observationforequation}$ denotes the Dirac-Delta functional associated to $\observationforequation$ and $\{\stdbase_i\}^{\Nobservations}_{i=1}$ is the canonical basis of $\R^{\Nobservations}$. Consequently, we define the linear and continuous \textit{parameter-to-state} mapping $\pts \colon D \to  C([T_0, T] \times \bar{\Omega}_o) $ as well as the \textit{parameter-to-observable} operator $\pto \colon D \to  \R^{\Nobservations} $ via $\pts(m)=u$ and $\pto(m)=\obsO \circ \pts(m)$, respectively. Given a \textit{misfit} vector $\misfit \in \mathbb{R}^{\Nobservations}$, e.g., $\misfit = \fracSolidus{1}{\sigma^2} (\obsO(u) - \measurement)$ with measurements $\measurement$ recorded under \tcr{white} noise, the associated \textit{misfit-to-adjoint} map is given by $\mta(\misfit) = \adjoint$, where $\adjoint$ is the solution to \tcr{the final value problem}
\begin{equation}\label{eq:adjoint_equation}
\begin{aligned}
  -\adjoint_t - \kappa \Delta \adjoint - \operatorname{div}(\adjoint \velocity) &= \sum_{i=1}^{\Nobservations} \misfiti[i] \, \delta_{\observationforequation} 
  && \text{in } (0, T) \times \Omega, \\
  (\velocity \adjoint + \kappa \nabla \adjoint) \cdot \normal &= 0 
  && \text{on } (0, T) \times (\outflowBoundary \cup \characteristicBoundary), \\
  \adjoint &= 0 
  && \text{on } (0, T) \times \inflowBoundary, \\
  \adjoint(T, \cdot) &= 0 
  && \text{in } \Omega. 
\end{aligned}
\tag{$\mathcal{P}_{\mta}$}
\end{equation}
We illustrate it in \autoref{fig:adjoint}.

\subsection{A sparse inversion framework for contaminant release} \label{sec:sparseReg}
For the specific context of airborne contaminant transport, it is reasonable to assume that sources are confined to a finite number of small, localized areas. In the present work, we incorporate this prior knowledge by modeling the sought-for sources, the \tcr{ground-truth}, as 
\begin{equation*}
\parameterInitial^\dagger=\sum^{N^\dagger_{\text{I}}}_{i=1}\,\lambda^{\text{I}, \dagger}_i\,\ansatzSourcesInitial(\x^{\text{I},\dagger}_i,\cdot), \quad \parameterContinuous^\dagger=\sum^{N^\dagger_{\text{C}}}_{j=1}\,\lambda^{\text{C}, \dagger}_j\,\ansatzSourcesContinuous(\x^{\text{C},\dagger}_j,\cdot), \quad \text{where} \quad (\xbr^{\text{I},\dagger},\sourceIntensity^{\text{I},\dagger}) \in \bar{\Omega}^{N^\dagger_{\text{I}}} \times \R^{N^\dagger_{\text{I}}}_{>0 }, \, (\xbr^{\text{C},\dagger},\sourceIntensity^{\text{C},\dagger}) \in \bar{\Omega}^{N^\dagger_{\text{C}}} \times \R^{N^\dagger_{\text{C}}}_{>0 }
\end{equation*}
denote unknown location-intensity pairs while $\ansatzSourcesInitial$ as well as $\ansatzSourcesContinuous$ are given shape functions, respectively. For a discussion of examples relevant to the present paper, we refer to \autoref{sec:shapeFunctions}. \tcr{Here, and with a slight abuse of notation, we interpret elements in $\bar{\Omega}^{N}$, $N\in \N$, as matrices whose columns are points in $\bar \Omega$.}

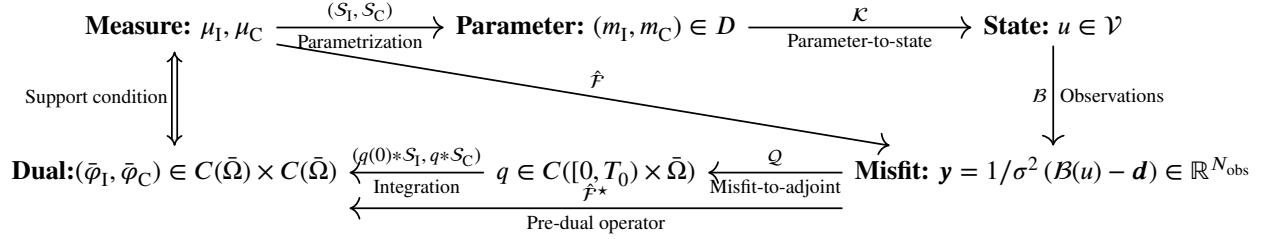
\begin{figure}
\makeatletter
\long\def\@makecaption#1#2{%
  \vskip\abovecaptionskip
  \noindent\textbf{#1:}\hspace{0.5em}#2\par
  \vskip\belowcaptionskip
}
\makeatother
\begin{tikzcd}[row sep=huge, column sep=large]
\textbf{Measure: } \measureInitial,\measureContinuous
    \arrow[d, Leftrightarrow, "\text{Support condition}" left]
    \arrow[r, "(\ansatzSourcesInitial{,}\,\ansatzSourcesContinuous)" above, "\text{Parametrization}" below] 
    \arrow[drr, "\hat \pto"]
& \textbf{Parameter: } (\parameterInitial,\parameterContinuous) \in\sop 
    \arrow[r, "\pts" above, "\text{Parameter-to-state}" below] 
& \textbf{State: } \concentration \in \bochner 
    \arrow[d, "\obsO" left, "\text{Observations}" right] \\
\textbf{Dual:} (\bar{\varphi}_{\text{I}}, \bar{\varphi}_{\text{C}}) \in C(\bar{\Omega}) \times C(\bar{\Omega})  
    \arrow[r,leftarrow, "(\adjoint(0)*\ansatzSourcesInitial{,}\,\adjoint*\ansatzSourcesContinuous)" above ,"\text{Integration}" below] 
\arrow[rr,leftarrow, shift right=3ex, "\hat \pto^\star" above ,"\text{Pre-dual operator}" below]
& \adjoint \in  C([0, T_0) \times \bar \Omega)
    \arrow[r,leftarrow,"\mta" above ,"\text{Misfit-to-adjoint}" below]
& \textbf{Misfit: } \misfit =\fracSolidus{1}{\sigma^2} (\obsO(u) - \measurement) \in \mathbb{R}^{\Nobservations}
\end{tikzcd}
\caption{\tcr{Summary of the inverse problem with sparse regularization}}
\label{fig:summary_measure}
\end{figure}

At this point, we emphasize that we neither assume knowledge of the positions and amplitudes nor of the number of sources. In contrast to the setting of \autoref{subsec:Ltworeg}, this leads to a highly nonlinear parameter identification problem. We alleviate these difficulties, by first noting that
\begin{equation}
\label{eq:finiteSum}
\parameterInitial^\dagger= \int_{\bar{\Omega}} \ansatzSourcesInitial(\x, \cdot)~\mathrm{d}\measureInitial^\dagger(\x), \quad \parameterContinuous^\dagger= \int_{\bar{\Omega}} \ansatzSourcesInitial(\x, \cdot)~\mathrm{d}\measureContinuous^\dagger (\x), \quad \text{where} \quad \measureInitial^\dagger =\sum^{N^\dagger_{\text{I}}}_{i=1}\,\lambda^{\text{I},\dagger}_i\,\delta_{\x^{\text{I},\dagger}_i},\,\measureContinuous^\dagger =\sum^{N^\dagger_{\text{C}}}_{j=1}\,\lambda^{\text{C}, \dagger}_j\,\delta_{\x^{\text{C},\dagger}_j}
\end{equation}
are conic combinations of Dirac-Delta functionals. Subsequently, we propose to recover the representing measures $(\measureInitial^\dagger, \measureContinuous^\dagger)$ rather than the sources $(\parameterInitial^\dagger, \parameterContinuous^\dagger)$ by solving
\begin{equation}\label{eq:sparseObjective}
\min_{\measureInitial,\measureContinuous \in \som(\bar\Omega)} J(\measureInitial,\measureContinuous) := \left \lbrack \fracSolidus{1}{(2\,\sigma^2)}\,\| \pto(\parameterInitial(\measureInitial),\parameterContinuous(\measureContinuous)) - \measurementPlain \|^2_{\R^{\Nobservations}} + \alpha\,\left(\measureInitial(\bar{\Omega})+\measureContinuous(\bar{\Omega})\right) \right \rbrack, \tag{$\mathcal{P}_{\mathcal{M}}$}
\end{equation} 
where $\alpha>0$ is a regularization parameter, $\sigma$ characterizes the measurement noise and
\begin{equation} \label{eq:linearparammaps}
\parameterInitial{(\measureInitial)}= \int_{\bar{\Omega}} \ansatzSourcesInitial(\x, \cdot)~\mathrm{d}\measureInitial(\x), \quad \parameterContinuous{(\measureContinuous)}= \int_{\bar{\Omega}} \ansatzSourcesContinuous(\x, \cdot)~\mathrm{d}\measureContinuous(\x).
\end{equation}
Here, $\som(\bar{\Omega})$ denotes the convex cone of positive Radon measures on $\bar{\Omega}$ and $\mu(\bar{\Omega})$ is the total variation of $\mu \in \som(\bar{\Omega})$. On the one hand, we emphasize that the mappings $\parameterInitial(\cdot)$ and $\parameterContinuous(\cdot)$ defined in \eqref{eq:linearparammaps} are linear. Consequently, the optimization problem \eqref{eq:sparseObjective} is a convex minimization problem and the existence of minimizers can be proven via the direct method in the calculus of variations. On the other hand, sparse measures, i.e. arbitrary conic combinations of Dirac-Delta functionals,
\begin{equation*}
    \somN(\bar{\Omega}):= \left\{\,\mu\lbrack\xbr, \sourceIntensity\rbrack=\sum^N_{j=1} \sourceIntensityi[j]\, \delta_{x_j}\;|\;\xbr\in \bar{\Omega}^N,~\sourceIntensity \in \R_+^N,~N \in \N\,\right\} \tcr{\text{ for which } \mu\lbrack\xb, \sourceIntensity\rbrack(\Omega)=|\sourceIntensity|_{\ell^1}\coloneq\sum_{j=1}^N |\sourceIntensityi[j]|} 
\end{equation*}
holds, form a proper, but not suitably closed subset of $\som(\bar{\Omega})$. Nevertheless, the following theorem establishes the existence of sparse minimizers to \eqref{eq:sparseObjective}, \tcr{thereby bridging the gap between our convex optimization framework and the structural assumption on the ground truth.} Moreover, it provides a tool to verify optimality. For this purpose, we introduce the \textit{lifted-parameter-to-observable map} $\hat{\pto}$, defined by $\hat{\pto}(\measureInitial, \measureContinuous)=\pto(\parameterInitial(\measureInitial), \parameterContinuous(\measureContinuous)),$
along with its formal pre-dual $\hat{\pto}^\star \colon \R^{\Nobservations} \to C(\bar{\Omega}) \times C(\bar{\Omega})$. Given $\misfit \in \R^{\Nobservations} $, the latter satisfies
\tcr{\begin{equation} \label{eq:valueadj}
\misfit^\top \hat{\pto}(\measureInitial, \measureContinuous)= \int_{\bar{\Omega}} \varphi_{\text{I}}(x)~\mathrm{d}\measureInitial(x) +\int_{\bar{\Omega}} \varphi_{\text{C}}(x)~\mathrm{d}\measureContinuous(x) 
\end{equation}}
where $ (\varphi_{\text{I}},\varphi_{\text{C}})\coloneqq \hat{\pto}^\star \misfit$ satisfies
\begin{equation}\label{eq:adjoint_term}
\left\lbrack\hat{\pto}^\star \misfit\right\rbrack (\x) = \left(
\int_{\bar\Omega} \ansatzSourcesInitial(\x,z) \, \mta(\misfit)(0,z)~\mathrm{d}\Omega(z),
\int_0^T \int_{\bar\Omega} \ansatzSourcesContinuous(\x,z) \, \mta(\misfit)(t,z)~\mathrm{d}\Omega(z)~\mathrm{dt}
\right) \quad \text{for all} \quad x \in \bar{\Omega}.
\end{equation}
Here we recall that $\mta(\misfit)$ is given by the solution to an advection-diffusion problem with inverted wind direction $\velocity$ (cf. \eqref{eq:adjoint_equation}).
\begin{theorem}\label{thm:optimality}
Problem \eqref{eq:sparseObjective} admits at least one minimizing pair $(\bmeasureInitial, \bmeasureContinuous) \in \somN(\bar \Omega)^2$, i.e.,
\begin{equation} \label{eq:sparsesolution}
    \bmeasureInitial= \measureInitial \lbrack\bar{\xbr}^{\text{I}},\bar{\sourceIntensity}^{\text{I}} \rbrack,\quad   \bmeasureContinuous=\measureContinuous \lbrack \bar{\xbr}^{\text{C}},\bar{\sourceIntensity}^{\text{C}}\rbrack, \quad \text{where} \quad  (\bar{\xbr}^{\text{I}},\bar{\sourceIntensity}^{\text{I}}) \in \bar{\Omega}^{\bar{N}_{\text{I}}} \times \Rplus^{\bar{N}_{\text{I}}}, \quad (\bar{\xbr}^{\text{C}},\bar{\sourceIntensity}^{\text{C}}) \in \bar{\Omega}^{\bar{N}_{\text{C}}} \times \Rplus^{\bar{ N}_{\text{C}}}, \quad \bar{N}_{\text{I}},\bar{N}_{\text{C}} \geq 0
\end{equation}
as well as $\bar{N}_{\text{I}}+\bar{N}_{\text{C}} \leq \Nobservations$.
Moreover, the following equivalence holds:   
\begin{itemize}
    \item A pair $(\bmeasureInitial, \bmeasureContinuous)$ of the form \eqref{eq:sparsesolution} is a minimizer of Problem \eqref{eq:sparseObjective}.   
    \item For the dual variables
\end{itemize} 
    \begin{equation} \label{eq:defdual}
        (\bar{\varphi}_{\text{I}}, \bar{\varphi}_{\text{C}}) \coloneqq - \hat{\pto}^\star \bar{{\misfit}} \in C(\bar{\Omega})\times C(\bar{\Omega}), \quad \text{where}\quad \bar{{\misfit}}=\fracSolidus{1}{\sigma^2} (\hat \pto(\bmeasureInitial,\bmeasureContinuous) - \measurementPlain),
    \end{equation}
    the following conditions hold:
    \begin{equation} \label{eq:suppcondition}
        \max_{\x \in \bar{\Omega}} \bar{\varphi}_{\text{I}}(\x) \leq \alpha, \quad \max_{\x \in \bar{\Omega}} \bar{\varphi}_{\text{C}}(\x) \leq \alpha, \quad \bar{\varphi}_{\text{I}}(\bar \x^{\text{I}}_i) = \alpha, \quad \bar{\varphi}_{\text{C}}(\bar{\x}^{\text{C}}_j) = \alpha \quad \text{for all} \quad  i=1, \dots \bar{N}_{\text{I}}, \quad j=1, \dots \bar{N}_{\text{C}}.
    \end{equation}
\end{theorem}
\begin{proof}
We start by discussing the claimed equivalence. For this purpose, set
\begin{equation*}
    f(\measureInitial,\measureContinuous) \coloneqq \fracSolidus{1}{(2\,\sigma^2)}\,\| \pto(\parameterInitial(\measureInitial),\parameterContinuous(\measureContinuous)) - \measurementPlain \|^2_{\R^{\Nobservations}}
\end{equation*}
By means of the chain rule as well as \eqref{eq:valueadj}, $f$ is Gateaux-differentiable at every pair $(\measureInitial, \measureContinuous) \in \som(\Omega)^2$ and there holds
\begin{equation*}
    f'(\measureInitial, \measureContinuous)(\delta\measureInitial, \delta\measureContinuous)= -\left(\int_{\bar{\Omega}} \varphi_{\text{I}}(x)~\mathrm{d}\delta\measureInitial(x) +\int_{\bar{\Omega}} \varphi_{\text{C}}(x)~\mathrm{d}\delta\measureContinuous(x) \right),
\end{equation*}
where $(\delta\measureInitial, \delta\measureContinuous)$ are two (signed) Radon measures and where we set
    \begin{equation} \label{eq:genericduals}
        ({\varphi}_{\text{I}}, {\varphi}_{\text{C}}) \coloneqq - \hat{\pto}^\star {\tcr{\misfit}} \in C(\bar{\Omega})\times C(\bar{\Omega}), \quad \text{with}\quad {\tcr{\misfit}}=\fracSolidus{1}{\sigma^2} (\hat \pto(\measureInitial,\measureContinuous) - \measurementPlain).
    \end{equation}
Since \eqref{eq:sparseObjective} is a convex problem, an arbitrary pair $(\bmeasureInitial, \bmeasureContinuous) \in \som(\Omega)^2$ is a minimizer if and only if the associated dual variables, \eqref{eq:defdual}, satisfy
\begin{equation*}
        \max_{\x \in \bar{\Omega}} \bar{\varphi}_{\text{I}}(\x) \leq \alpha, \quad \max_{\x \in \bar{\Omega}} \bar{\varphi}_{\text{C}}(\x) \leq \alpha, \quad \operatorname{supp} \bmeasureInitial \subset \left\{\,x\in \bar\Omega\;|\;\bar{\varphi}_{\text{I}}(x)=\alpha\,\right\}, \quad \operatorname{supp} \bmeasureContinuous \subset \left\{\,x\in \bar \Omega\;|\;\bar{\varphi}_{\text{C}}(x)=\alpha\,\right\}
    \end{equation*}
where $\operatorname{supp} \bmeasureInitial$, $\operatorname{supp} \measureContinuous$ denote the support of $\bmeasureInitial$ and $\bmeasureContinuous$, respectively. This can be derived, e.g., analogously to \cite[Proposition 3.5,~Lemma 3.6]{neitzel19}.
We conclude the claimed equivalence noting that for measures $(\bmeasureInitial, \bmeasureContinuous)$ of the form \eqref{eq:sparsesolution}, there holds
\begin{equation*}
    \operatorname{supp} \bmeasureInitial= \left\{\,\bar{x}^{\text{I}}_i\;|\;i=1,\dots, \bar{N}_{\text{I}}\,\right\}, \quad  \operatorname{supp} \bmeasureContinuous= \left\{\,\bar{x}^{\text{C}}_j\;|\;j=1,\dots, \bar{N}_{\text{C}}\,\right\}.
\end{equation*}
It thus remains to show that sparse minimizers with $\bar{N}_{\text{I}}+\bar{N}_{\text{C}} \leq \Nobservations$ exist. By applying an abstract convex representer theorem, see e.g. \cite{Boyer.2019, Carioni.2020}, 
\eqref{eq:sparseObjective} admits at least one solution of the form
\begin{equation*}
    (\bmeasureInitial, \bmeasureContinuous)= \sum^{\bar N}_{j=1} \bar \lambda_j (\measureInitial^j, \measureContinuous^j), \quad \text{with} \quad \bar{\lambda}_j \geq 0,~\bar{N} \leq \Nobservations \quad \text{and} \quad (\measureInitial^j, \measureContinuous^j) \in \operatorname{Ext}(B).
\end{equation*}
Here $\operatorname{Ext}(B)$ denotes the extremal points of the weak*-compact set
\begin{equation*}
    B \coloneqq \left\{\,(\mu_{\text{I}}, \mu_{\text{C}}) \in \som(\bar\Omega)^2\;|\; \alpha \left( \mu_{\text{I}}(\bar \Omega)+\mu_{\text{C}}(\bar \Omega)\right) \leq 1\,\right\}.
\end{equation*}
Following arguments similar to \cite[Lemma 3.2]{Bredies.2025}, we obtain
\begin{equation*}
\operatorname{Ext}(B)= \{(0,0)\} \cup \left\{\,\frac{1}{\alpha}\left(\delta_{x}, 0 \right)\;|\;x \in \bar\Omega\,\right\}\cup \left\{\,\frac{1}{\alpha}\left(0,\delta_{x} \right)\;|\;x \in \bar \Omega\,\right\},
\end{equation*}
which proves the sparsity of the minimizer, see (Eq. \ref{eq:sparsesolution}).
\end{proof}
\begin{remark}
    We emphasize that this is an existence but not a uniqueness result. Quite the contrary, since $\pto$ has finite rank, Problem \eqref{eq:sparseObjective} can admit a multitude of solutions, some of which might even be non-sparse, i.e. elements of $\som(\bar \Omega ) \setminus \somN(\bar \Omega)$. Additionally, identifying the correct number of sources, $\bar N_{\text{I}}= N^\dagger_{\text{I}}$ and $\bar N_{\text{C}}= N^\dagger_{\text{C}}$, together with quantitative reconstruction results requires strong nondegeneracy conditions, small measurement noise as well as a noise-adapted regularization parameter $\alpha$. We refer, e.g., to \cite{Duval.2015} where these topics are addressed in detail for simpler deconvolution tasks.
\end{remark}
We close this section by noting that the support condition in \eqref{eq:suppcondition} can be used to verify optimality of sparse candidate measures. Moreover, for sparse minimizers, the objective functional decomposes as 
\begin{equation*}
    J(\bmeasureInitial, \bmeasureContinuous)=J(\measureInitial \lbrack\bar{\xbr}^{\text{I}},\bar{\sourceIntensity}^{\text{I}}\rbrack, \measureContinuous \lbrack \bar{\xbr}^{\text{C}},\bar{\sourceIntensity}^{\text{C}}\rbrack)= \fracSolidus{1}{(2\,\sigma^2)}\,\|\hat{\pto}(\measureInitial \lbrack\bar{\xbr}^{\text{I}},\bar{\sourceIntensity}^{\text{I}}\rbrack, \measureContinuous \lbrack \bar{\xbr}^{\text{C}},\bar{\sourceIntensity}^{\text{C}}\rbrack) - \measurementPlain\|^2_{\R^{\Nobservations}} + \alpha\,\left(|\bar{\sourceIntensity}^{\text{I}}|_{\ell_1}+|\bar{\sourceIntensity}^{\text{C}}|_{\ell_1})\right),
\end{equation*}
where
$\hat{\pto}(\measureInitial \lbrack\bar{\xbr}^{\text{I}},\bar{\sourceIntensity}^{\text{I}}\rbrack,\measureContinuous \lbrack \bar{\xbr}^{\text{C}},\bar{\sourceIntensity}^{\text{C}}\rbrack)= \sum^{\bar{N}_{\text{I}}}_{i=1}  \bar{\lambda}^{\text{I}}_i \hat{\pto}\left(\delta_{\bar{\x}^\text{I}_i},0\right) + \sum^{\bar{N}_{\text{C}}}_{j=1} \bar{\lambda}^{\text{C}}_j \hat{\pto}\left(0,\delta_{\bar{\x}^\text{C}_j}\right)$
is linear w.r.t to the source intensities. Hence, once optimal locations have been identified, the associated intensities are obtained as minimizers of a finite-dimensional, convex problem.
\begin{corollary} \label{coroll:finiteopt}
Consider a minimizing pair $(\bmeasureInitial, \bmeasureContinuous)$ of the form \eqref{eq:sparsesolution}. Then there holds
\begin{equation} \label{eq:charabyfinite}
(\bar{\sourceIntensity}^{\text{I}},\bar{\sourceIntensity}^{\text{C}}) \in \argmin_{\sourceIntensity^\text{I} \in \R^{N_{\text{I}}}_{\geq 0},\sourceIntensity^\text{C} \in \R^{{N}_{\text{C}}}_{\geq 0} } \left\lbrack \fracSolidus{1}{(2\,\sigma^2)}\,\|\hat{\pto}(\measureInitial \lbrack\bar{\xbr}^{\text{I}},{\sourceIntensity}^{\text{I}}\rbrack, \measureContinuous \lbrack \bar{\xbr}^{\text{C}},{\sourceIntensity}^{\text{C}}\rbrack) - \measurementPlain\|^2_{\R^{\Nobservations}} + \alpha\,\left(|{\sourceIntensity}^{\text{I}}|_{\ell_1}+|{\sourceIntensity}^{\text{C}}|_{\ell_1})\right) \right\rbrack.
\end{equation}
Vice versa, let two arbitrary sparse measures of the form \eqref{eq:sparsesolution} be given, define $(\bar{\varphi}_{\text{I}}, \bar{\varphi}_{\text{C}})$ according to \eqref{eq:defdual} and assume that the pairs $(\bar{\xbr}^{\text{I}},\bar{\sourceIntensity}^{\text{I}}) \in \bar{\Omega}^{\bar{N}_{\text{I}}} \times \Rplus^{\bar{N}}, (\bar{\xbr}^{\text{C}},\bar{\sourceIntensity}^{\text{C}}) \in \bar{\Omega}^{\bar{N}_{\text{C}}} \times \Rplus^{\bar{ N}_{\text{C}}}$
satisfy \eqref{eq:charabyfinite}. Then we have
\begin{equation*}
    (\bmeasureInitial, \bmeasureContinuous)~\text{is a minimizer of}~\eqref{eq:sparseObjective} \Leftrightarrow \max_{\x \in \bar{\Omega}} \bar{\varphi}_{\text{I}}(\x) \leq \alpha \text{ and } \max_{\x \in \bar{\Omega}} \bar{\varphi}_{\text{C}}(\x) \leq \alpha.
\end{equation*}
\end{corollary}
\tcr{\begin{proof}
With a slight abuse of notation, define the mapping
\begin{equation*}
f(\sourceIntensity^{\text{I}},\sourceIntensity^{\text{C}}) \coloneqq  \fracSolidus{1}{(2\,\sigma^2)}\,\|\hat{\pto}(\measureInitial \lbrack\bar{\xbr}^{\text{I}},{\sourceIntensity}^{\text{I}}\rbrack, \measureContinuous \lbrack \bar{\xbr}^{\text{C}},{\sourceIntensity}^{\text{C}}\rbrack) - \measurementPlain\|^2_{\R^{\Nobservations}} \
\end{equation*}
and note that
\begin{equation*}
    J(\measureInitial \lbrack\bar{\xbr}^{\text{I}},{\sourceIntensity}^{\text{I}}\rbrack, \measureContinuous \lbrack \bar{\xbr}^{\text{C}},{\sourceIntensity}^{\text{C}}\rbrack)=f(\sourceIntensity^{\text{I}},\sourceIntensity^{\text{C}})+\alpha\,\left(|{\sourceIntensity}^{\text{I}}|_{\ell_1}+|{\sourceIntensity}^{\text{C}}|_{\ell_1}\right).
\end{equation*}
For a minimizing pair $(\bmeasureInitial, \bmeasureContinuous)$ of the form \eqref{eq:sparsesolution}, we now estimate
\begin{equation*}
    J(\measureInitial \lbrack\bar{\xbr}^{\text{I}},\bar{\sourceIntensity}^{\text{I}}\rbrack, \measureContinuous \lbrack \bar{\xbr}^{\text{C}},\bar{\sourceIntensity}^{\text{C}}\rbrack)=    J(\bmeasureInitial, \bmeasureContinuous) \leq J(\measureInitial \lbrack\bar{\xbr}^{\text{I}},{\sourceIntensity}^{\text{I}}\rbrack, \measureContinuous \lbrack \bar{\xbr}^{\text{C}},{\sourceIntensity}^{\text{C}}\rbrack) \quad \text{for all} \quad \sourceIntensity^\text{I} \in \R^{N_{\text{I}}}_{\geq 0},\sourceIntensity^\text{C} \in \R^{{N}_{\text{C}}}_{\geq 0}
\end{equation*}
from which we conclude \eqref{eq:charabyfinite}.
For the second statement, and again by means of the chain rule, note that $f$ is Gateaux-differentiable at every pair $(\sourceIntensity^\text{I},\sourceIntensity^\text{C}) \in  \R^{ \bar N_{\text{I}}}_{\geq 0} \times \R^{{N}_{\text{C}}}_{\geq 0}$ and there holds
\begin{equation*}
f'(\sourceIntensity^\text{I},\sourceIntensity^\text{C})(\delta \sourceIntensity^\text{I}, \delta\sourceIntensity^\text{C})= -\left( \sum^{\bar{N}_{\text{I}}}_{i=1} \delta \lambda^{\text{I}}_i \,\varphi_{\text{I}}(\bar{x}^{\text{I}}_i)+\sum^{\bar{N}_{\text{C}}}_{j=1} \delta \lambda^{\text{C}}_j\,\varphi_{\text{C}}(\bar{x}^{\text{C}}_J)  \right )
\end{equation*}
for all directions $(\delta \sourceIntensity^\text{I}, \delta\sourceIntensity^\text{C})  \in \R^{ \bar N_{\text{I}}} \times \R^{{N}_{\text{C}}}$ and where we define
    \begin{equation*}
        ({\varphi}_{\text{I}}, {\varphi}_{\text{C}}) \coloneqq - \hat{\pto}^\star {\tcr{\misfit}} \in C(\bar{\Omega})\times C(\bar{\Omega}), \quad \text{with}\quad {\tcr{\misfit}}=\fracSolidus{1}{\sigma^2} (\hat \pto(\measureInitial \lbrack\bar{\xbr}^{\text{I}},{\sourceIntensity}^{\text{I}}\rbrack, \measureContinuous \lbrack \bar{\xbr}^{\text{C}},{\sourceIntensity}^{\text{C}}\rbrack) - \measurementPlain).
    \end{equation*}
    analogously to \autoref{thm:optimality}. Since \eqref{eq:charabyfinite} is also a convex minimization problem, optimality of $(\bar{\sourceIntensity}^{\text{I}},\bar{\sourceIntensity}^{\text{C}})$ is equivalent to
    \begin{equation*}
        \bar{\varphi}_{\text{I}}(\bar \x^{\text{I}}_i) = \alpha, \quad \bar{\varphi}_{\text{C}}(\bar{\x}^{\text{C}}_j) = \alpha \quad \text{for all} \quad  i=1, \dots \bar{N}_{\text{I}}, \quad j=1, \dots \bar{N}_{\text{C}}.
    \end{equation*}
 which follows by standard first-order theory, see e.g. \cite{boyd}. Comparing with the second part of \eqref{thm:optimality}, the claimed equivalence follows.
\end{proof}}
The connection between the different components of the sparse regularization approach is illustrated in \autoref{fig:summary_measure}.
\subsection{{Modeling} of contaminant sources}
\label{sec:shapeFunctions}
\begin{figure}
\centering
\begin{subfigure}{.5\textwidth}
\includegraphics[width=1.0\linewidth]{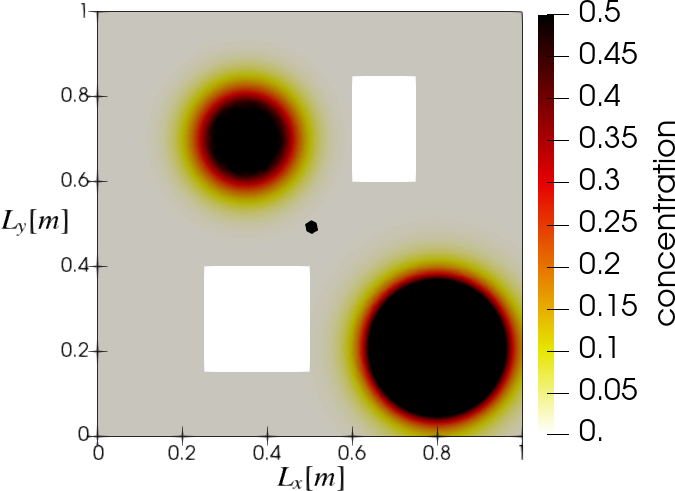}
\end{subfigure}
\caption{Mathematical description of contaminant release. Closed form expression $\ansatzSources^1([0.35, 0.7],0.26,\cdot)$, elliptic PDE solution $\ansatzSources^2([0.8, 0.2],1.0,0.001,\cdot)$, and very localized Dirac distribution ${\ansatzSources^3([0.5, 0.5],\cdot)}$.}
\label{fig:initial_types}
\end{figure}

In this paper, three models of contaminant release are presented and the type of the problem-defining parameter $m = (\parameterInitial, \parameterContinuous)$ is specified accordingly. Numerical examples of the three types are illustrated in \autoref{fig:initial_types}. The first model is a radial basis function given in closed form expression as
\begin{equation}
 \begin{aligned}
  \ansatzSources^1(\x_s,r,y)=\min \left\{0.5,\exp \left(-\text{ln}(\epsilon)\norm{y-\x_s}_2^2/r^2\right)\right\}.
\end{aligned}\label{eq:gauss_blob}
\end{equation}
Here $\x_s \in \Omega$ is the center, $r>0$ is the radius of the source and $\epsilon>0$ is a given threshold. The model $\ansatzSources^1 $ is employed for inverse problems with $L^2$-regularization \cite{Villa.2021,Petra.2011} and is a highly flexible approach that can be tailored to the specific requirements of the problem. The second model $\ansatzSources^2$ for a contaminant release is constructed as solution of an elliptic PDE with a Dirac distribution with center $\x_s$ on the right-hand side:
\begin{equation}
\begin{aligned}
\ansatzSources^2(\x_s, \ellipI, \ellipLaplace,y) = m_{\x_s}(y), 
\text{ with }m_{\x_s} \text{ given by } \quad \left\{ \begin{array}{llll}
(\ellipI I - \ellipLaplace \Delta) m_{\x_s} &= \delta_{\x_s} &\quad &\text{in } \Omega, \\
\ellipLaplace \nabla m_{\x_s} \cdot \normal + \beta m_{\x_s} &= 0 & &\text{on } \partial \Omega,
\end{array}\right.
\end{aligned}
\label{eq:elliptic_initial}
\end{equation}
where the scalars $\ellipI$, $\ellipLaplace$ and  $\beta(\eta,\gamma)$ are used to adjust the initial condition to the specific application. The last formulation 
\begin{equation}
\begin{aligned}
  \ansatzSources^3(\x_s)= \delta_{x_s}
\end{aligned}
\label{eq:s3}
\end{equation}
is particularly advantageous for large-scale problems in which the source has minimal support.

\subsection{Quadratic regularization in source space}
\label{subsec:Ltworeg}
A natural alternative to the parametrization ansatz outlined above is to consider a regularization term directly acting on the source space $D$ leading to problems of the form
    \begin{equation}\label{eq:regObjective}
\min_{(\parameterInitial, \parameterContinuous) \in D} \left \lbrack \fracSolidus{1}{(2\,\sigma^2)}\,\| \pto(\parameterInitial,\parameterContinuous) - \measurement\|^2_{\R^{\Nobservations}} + \mathcal{R}(m) \right \rbrack. \tag{$\mathcal{P}_D$}
\end{equation}
In this regard, a popular approach is to consider $D \subset L^2(\Omega) \times L^2(\Omega) $ and choose quadratic regularization terms of the form $\mathcal{R}(m):=\fracSolidus{1}{2}\norm{\mathcal{A}(m-m_{\text{prior}})}^2_{L^2(\Omega)^2}$ for a reference $m_{\text{prior}}$ and, e.g., a differential operator $\mathcal{A}$. This admits a Bayesian interpretation: Assuming that $m \sim \mathcal{N}(m_{\text{prior}}, \mathcal{A}^{-2}) $ with $\mathcal{A}^{-2}$ a trace-class operator and $m_{\text{prior}}$ sufficiently regular as well as a linear measurement model $d= \pto(m)+\noise$, $\noise \sim \mathcal{N}(0, \sigma^2 \operatorname{Id}$), \eqref{eq:regObjective} corresponds to the computation of the maximum a posteriori estimator. Since this problem is convex, the associated necessary and sufficient first-order optimality condition implies $\bar m = \mathcal{H}^{-1}\left((\fracSolidus{1}{\sigma^2})\pto\,\measurement + \mathcal{A}^{2}m_{\text{prior}}\right)$ where $\mathcal{H}=(\fracSolidus{1}{\sigma^2})\pto^*\pto + \mathcal{A}^{2}$. For a given realization of the measurement data, an estimate can thus be obtained by solving this linear equation using an inexact Newton CG method~\cite{Steihaug.1983} for a suitable discretization~\cite{Villa.2021,Petra.2011}. Due to the trace-class assumption, the eigenvalues of the positive-definite matrix $\mathcal{H}$ decay rapidly enabling the use of reduced models. This property makes real-time parameter identification feasible through an offline-online decomposition strategy. Restricting ourselves to the case of initial value identification, i.e., assuming $\parameterContinuous = 0$, we choose a Laplacian-like operator ${\mathcal{A} := (\ellipI\, I - \ellipLaplace \Delta)}$, cf.~\cite{Villa.2021,Petra.2011}, together with Robin boundary conditions $\ellipLaplace \nabla \parameterInitial \cdot \normal + \beta \parameterInitial = 0 \text{ in } (0,T) \times \partial \Omega,$ with constants as in \cite{Daon.2018}. An available implementation of the described method, see \cite{Wu.2023}, is employed as a benchmark for the computations in \autoref{sec:numerical}.

\section{Primal-Dual-Active-Point Algorithm for Source Identification} \label{sec:pdap}
\begin{figure}
\centering
\begin{subfigure}{0.48\textwidth}
\includegraphics[width=0.9\linewidth]{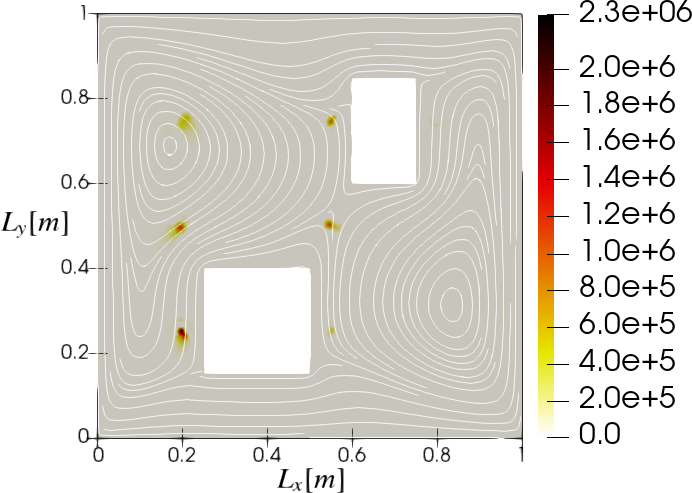}
\caption{Sensor misfit $\misfit$}
\end{subfigure}
\centering
\begin{subfigure}{0.48\textwidth}
\includegraphics[width=0.9\linewidth]{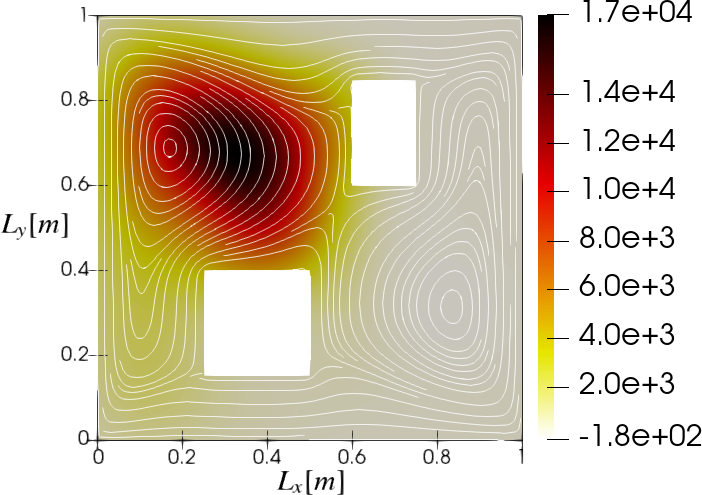}
\caption{Dual variable $\varphi_{\text{I}}$}
\end{subfigure}
\centering
\begin{subfigure}{0.48\textwidth}
\includegraphics[width=0.9\linewidth]{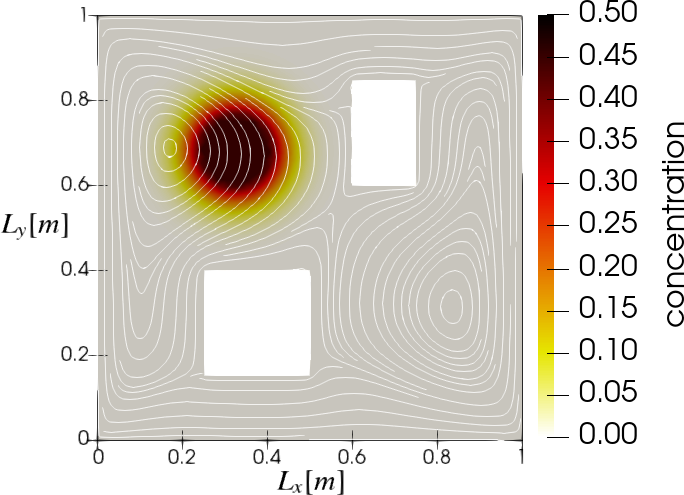}
\caption{Estimated initial condition $\parameterInitial$}
\end{subfigure}
\centering
\begin{subfigure}{0.48\textwidth}
\includegraphics[width=0.9\linewidth]{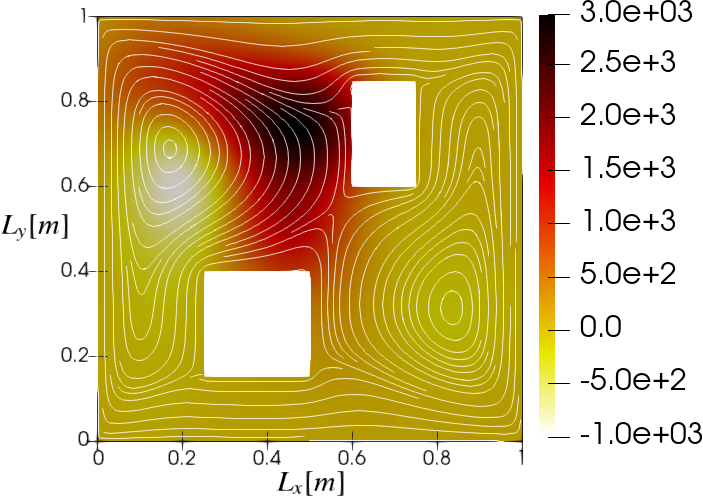}
\caption{Updated dual variable $\varphi_{\text{I}}$ }
\end{subfigure}
\caption{Iteration of PDAP-algorithm. Starting with sensor misfit $\misfit$ (a) obtained from ($\parameterInitial=0$), the dual variable $\varphi_{\text{I}}$ is computed (b). Position of $\max({\varphi_{\text{I}}})$ provides a candidate to estimate the initial condition $\parameterInitial$ after the first iteration (c). Reduced value of updated field $\varphi_{\text{I}}$ (d) indicates convergence towards stopping criterion $\max(\varphi_{\text{I}}) \leq \alpha$.}
\label{fig:adjointandinitial}
\end{figure}

Finally, we briefly discuss the numerical solution of \eqref{eq:sparseObjective} by adapting the Primal-Dual-Active-Point (PDAP) method as described in \cite{Pieper.2021}. In this regard, and in accordance with the original work, we follow an optimize-then-discretize philosophy and formulate the algorithm on the function space level before applying it to suitable discretizations described in \autoref{sec:discretization}. Conceptually, PDAP is a greedy algorithm which \tcr{relies on sparse iterates} and exploits the existence of a pair of sparse minimizers, see  \autoref{thm:optimality}, \tcr{by alternating between proposing new candidate sources based on the violation of the optimality condition \eqref{eq:suppcondition} and removing unnecessary points motivated by the observations in \autoref{coroll:finiteopt}}. More in detail, it constructs a sequence of parametrized iterates
\begin{equation*}
        \measureInitial^k= \measureInitial \lbrack{\xbr}^{\text{I}}_k,{\sourceIntensity}^{\text{I}}_k \rbrack,\quad   \measureContinuous^k=\measureContinuous \lbrack {\xbr}^{\text{C}}_k,{\sourceIntensity}^{\text{C}}_k\rbrack \quad \text{where} \quad  ({\xbr}^{\text{I}}_k,{\sourceIntensity}^{\text{I}}_k) \in \bar{\Omega}^{{N}^k_{\text{I}}} \times \Rplus^{{N}^k_{\text{I}}}, \quad ({\xbr}^{\text{C}}_k,{\sourceIntensity}^{\text{C}}_k) \in \bar{\Omega}^{{N}^k_{\text{C}}} \times \Rplus^{{N}^ k_{\text{C}}}
\end{equation*}
with the property that $({\sourceIntensity}^{\text{I}}_k,{\sourceIntensity}^{\text{C}}_k)$ is a minimizing pair of
\begin{equation} \label{eq:finitealgo}
 \min_{\sourceIntensity^\text{I} \in \mathbb{R}_{\geq 0}^{{N}^k_{\text{I}}}, \, \sourceIntensity^\text{C} \in \mathbb{R}_{\geq 0}^{{N}^k_{\text{C}}} } \left\lbrack \fracSolidus{1}{(2\,\sigma^2)}\,\|\hat{\pto}(\measureInitial \lbrack{\xbr}^{\text{I}},{\sourceIntensity}^{\text{I}}\rbrack, \measureContinuous \lbrack {\xbr}^{\text{C}},{\sourceIntensity}^{\text{C}}\rbrack) - \measurement\|^2_{\R^{\Nobservations}} + \alpha\,\left(|{\sourceIntensity}^{\text{I}}|_{\ell_1}+|{\sourceIntensity}^{\text{C}}|_{\ell_1})\right) \right\rbrack \tag{$\mathcal{P}(\xbr^{\text{I}},\xbr^{\text{C}})$}
\end{equation}
for the particular choice of $(\xbr^{\text{I}},\xbr^{\text{C}})=(\xbr^{\text{I}}_k,\xbr^{\text{C}}_k)$. Referring to the terminology of the original paper, cf.~Algorithm~1 in~\cite{Pieper.2021}, $\xbr^{\text{I}}_k$ and $\xbr^{\text{C}}_k$, respectively, collect the currently ``active'' points comprising the iterates. In each iteration, we start by computing the current misfit and dual variables
\begin{equation*}
  {\misfit_k}=  \fracSolidus{1}{\sigma^2} (\hat \pto( \measureInitial^k, \measureContinuous^k) - \measurementPlain) \in \mathbb{R}^{\Nobservations}, \quad  ({\varphi}^k_{\text{I}}, {\varphi}^k_{\text{C}}) = - \hat{\pto}^\star {\misfit_k} \in C(\bar{\Omega})\times C(\bar{\Omega})
\end{equation*}
determine global maximizers $\widehat{\xb}^I_k$ and $\widehat{\xb}^C_k$ of the latter and check for optimality according to \autoref{coroll:finiteopt}. If this is not achieved up to a tolerance, we append $\widehat{\xb}^I_k$ and/or $\widehat{\xb}^C_k$ to the location vectors yielding $(\xbr^{\text{I}}_{k+1/2},\xbr^{\text{C}}_{k+1/2})$, respectively. Subsequently, we adjust the iterate by optimizing the intensities 
\begin{equation*}
    \measureInitial^{k+1}= \measureInitial \left\lbrack\xbr^{\text{I}}_{k+1/2},{\sourceIntensity}^{\text{I}}_{k+1/2}\right\rbrack,\quad   \measureContinuous^{k+1}=\measureContinuous \left\lbrack \xbr^{\text{C}}_{k+1/2},{\sourceIntensity}^{\text{C}}_{k+1/2} \right\rbrack,\quad \left({\sourceIntensity}^{\text{I}}_{k+1/2},{\sourceIntensity}^{\text{C}}_{k+1/2}\right) \in \argmin \left(\text{\hyperref[eq:finitealgo]{{$\mathcal{P}\left(\xbr^{\text{I}}_{k+1/2},\xbr^{\text{C}}_{k+1/2}\right)$}}}\right).  
\end{equation*}
Finally, location and intensity vectors are pruned by removing entries corresponding to Dirac-Delta functionals with zero weights. We summarize this procedure in \autoref{alg:pdap}. An example for a typical iteration for the simplified case of initial value identification is displayed in \autoref{fig:adjointandinitial}.

\tcr{By adapting the results in \cite{Pieper.2021}, we point out that $(\measureInitial^k,\measureContinuous^k)$ weak*-converges, on subsequences, towards minimizers of Problem \eqref{eq:sparseObjective} with a global, sublinear worst-case rate of convergence, i.e. there holds
\begin{equation*}
  r_J(\measureInitial^k,\measureContinuous^k) \coloneqq J(\measureInitial^k,\measureContinuous^k) - \min_{\measureInitial,\measureContinuous \in \som(\Omega)} J(\measureInitial,\measureContinuous) \leq \fracSolidus{C}{k} \quad \text{for all} \quad k \geq 1
\end{equation*}
and some $C>0$. 
Additionally, the dual variables satisfy
\begin{equation*}
    0 \leq   \max_{\x \in \Omega} \varphi^k_{\text{I}}(\x) -\alpha\rightarrow 0, \quad 0 \leq \max_{\x \in \Omega} \varphi^k_{\text{C}}(\x) -\alpha \rightarrow 0.
\end{equation*}
Hence, for $\text{tol}>0$, \autoref{alg:pdap} converges in a finite number of steps $K$ and we have 
\begin{equation}\label{eq:residual}
r_J(\measureInitial^k,\measureContinuous^k)\leq \frac{J(0,0)}{\alpha} \text{tol}
\end{equation}
assuming that the subproblems in Step 3 and 6 of \autoref{alg:pdap} are solved to optimality, see, mutatis mutandis, \cite[Proposition 4.12]{walter.dissertation} and the discussion on \cite[p. 83]{walter.dissertation} . Moreover, asymptotic linear decrease, i.e. 
$
    r_J(\measureInitial^{k+1},\measureContinuous^{k+1}) \leq \zeta r_J(\measureInitial^k,\measureContinuous^k),
$
for all $k$ large enough and some $\zeta \in (0,1)$,
can be ensured under additional structural assumptions related to a suitable no-gap second-order condition, \cite{wachsmuthwalter}. The theoretical upper bound on the support size from \autoref{thm:optimality},  ${N}^k_{\text{I}}+{N}^k_{\text{C}} \leq \Nobservations$ can be ensured throughout all iterations by augmenting the method with additional sparsification strategies, see e.g., \cite[Algorithm 2]{neitzel19}. However, in practice, this is not required since the method usually terminates before this threshold is even reached.
}
\begin{algorithm}
\caption{Primal-Dual-Active-Point-Strategy for Source Identification}\label{alg:pdap}
\begin{algorithmic}
\Require $\ansatzSourcesInitial$,$\ansatzSourcesContinuous$,$\measureInitial^0=\measureContinuous^0=0$, \tcr{empty matrices and vectors}  $\xbr^{\text{I}}_k,{\sourceIntensity}^{\text{I}}_k,\xbr^{\text{C}}_k,{\sourceIntensity}^{\text{C}}_k$\\
\For{$k = 0,1,2...$}
\State~
\State 1.\,Given $\measureInitial^k= \measureInitial \lbrack{\xbr}^{\text{I}}_k,{\sourceIntensity}^{\text{I}}_k \rbrack,\,\measureContinuous^k=\measureContinuous \lbrack {\xbr}^{\text{C}}_k,{\sourceIntensity}^{\text{C}}_k\rbrack$, compute $
  {\misfit}_k =  \fracSolidus{1}{\sigma^2} \left(\hat \pto( \measureInitial^k, \measureContinuous^k) - \measurement\right).$
\State~
\State 2.\,Compute convolution
\begin{equation*}
({\varphi}^k_{\text{I}}, {\varphi}^k_{\text{C}}) =-\left(
\int_\Omega \ansatzSourcesInitial(\cdot,z) \, \mta({\misfit}_k)(0,z)~\mathrm{d}{\Omega}(z),
\int_0^T \int_\Omega \ansatzSourcesContinuous(\cdot,z) \, \mta({\misfit}_k)(t,z)~\mathrm{d}{\Omega}(z)~\mathrm{dt}\right). 
\end{equation*}
\State~
\State 3.\,Determine global maxima
\begin{equation*}
    \widehat{\x}_k^{\text{I}}\in \argmax_{\x \in \Omega} \varphi^k_{\text{I}}(\x), \quad \widehat{\x}_k^{\text{C}}\in \argmax_{\x \in \Omega} \varphi^k_{\text{C}}(\x).
\end{equation*}
\State~
\State 4.\, Set $\xbr^{\text{I}}_{k+1/2}=\xbr^{\text{I}}_{k}$,~$\xbr^{\text{C}}_{k+1/2}=\xbr^{\text{C}}_{k}$. 
\State~
\State 5.\, Append $\widehat{\x}_k^{\text{I}}$ to $\xbr^{\text{I}}_{k+1/2}$ \textbf{if} $ \varphi^k_{\text{I}}(\widehat{\x}_k^{\text{I}})> \alpha$ + \text{tol}, \quad append $\widehat{\x}_k^{\text{C}}$ to $\xbr^{\text{C}}_{k+1/2}$ \textbf{if} $ \varphi_{\text{C}}(\widehat{\x}_k^{\text{C}})> \alpha+ \text{tol}$
\State \,\,\quad 
\State 6.\, Update source intensities
\begin{equation*}
\left({\sourceIntensity}^{\text{I}}_{k+1/2},{\sourceIntensity}^{\text{C}}_{k+1/2}\right) \in \argmin \left(
\text{\hyperref[eq:finitealgo]{{$\mathcal{P}\left(\xbr^{\text{I}}_{k+1/2},\xbr^{\text{C}}_{k+1/2}\right)$}}}\right).
\end{equation*}
\State \,\,\quad
\State 7. Update iterates
\begin{equation*}
    \measureInitial^{k+1}= \measureInitial \left\lbrack\xbr^{\text{I}}_{k+1/2},{\sourceIntensity}^{\text{I}}_{k+1/2}\right\rbrack,\quad   \measureContinuous^{k+1}=\measureContinuous \left\lbrack \xbr^{\text{C}}_{k+1/2},{\sourceIntensity}^{\text{C}}_{k+1/2} \right\rbrack.
\end{equation*}
\State \,\,\quad
\State 8.\, Obtain $\xbr^{\text{I}}_{k+1},{\sourceIntensity}^{\text{I}}_{k+1}$ and $\xbr^{\text{C}}_{k+1},{\sourceIntensity}^{\text{C}}_{k+1}$ by pruning.
\State~
\EndFor
\end{algorithmic}
\end{algorithm}
\tcr{We briefly comment on the practical realization of \autoref{alg:pdap}.}
By saving
$\hat{\pto}\left(\delta_{\widehat{\x}^I_k},0\right)$ and $ \hat{\pto}\left(0,\delta_{\widehat{\x}^\text{C}_k}\right)$ across iterations and due to the linearity of $\pto$, every iteration requires at most two forward simulations as well as one solution of the adjoint equation. No additional PDE solves are required for the realization of \eqref{eq:finitealgo}. \tcr{Concerning the latter, a semismooth Newton method based on the normal-map approach with multidimensional filter globalization is implemented~\cite{Milzarek.2014}, which we warmstart using the current intensity vectors to construct an initial iterate in order to benefit from its potentially fast, local convergence. The Newton solver is stopped once an appropriate residual falls below $\text{tol}_{\text{Newt}}<< \text{tol}$.}

Next, obtaining the dual variables via integration against the shape functions might be expensive if $\ansatzSources$ has large support in $\Omega$. However, for the PDE-based release model $\ansatzSources^2$, given in Eq.~\eqref{eq:elliptic_initial}, the integration can be simplified via adjoint calculus leading to a single solve of an elliptic equation, while the Dirac-based formulation of contaminant release $\ansatzSources^3$, given in Eq.~\eqref{eq:s3}, renders this step trivial. \tcr{ Finally, we acknowledge that the point insertion steps, that is, the computation of global maxima in Step 3 of \autoref{alg:pdap}, are in general nonconcave maximization problems. Solving solutions that are (approximately) globally optimal is computationally infeasible for general $\Omega$. Approaches to address these challenges can be found in, e.g., \cite{Flinth25, Hnatiuk.2025}.  
In our case, as outlined above, we ultimately apply \autoref{alg:pdap} to discretizations of Problem \eqref{eq:sparseObjective}. This also entails replacing $ \som(\Omega)$ by its discretized counterpart, the cone spanned by Dirac-Deltas in the nodes of the finite elements mesh. This reduces Step 3 of \autoref{alg:pdap} to sorting the values of $\varphi_{\text{I}}$ and $\varphi_{\text{I}}$, respectively, in the grid nodes, see also \autoref{sec:discretization}.}
\section{Finite Element Discretization of Forward and Adjoint Problem} \label{sec:discretization}
\tcr{In this section, we derive the weak formulation of~\eqref{eq:forward_equation}. 
To this end, we introduce suitable spaces for test functions
$ H^{1,2}_{\inflowBoundary}(\Omega)
:= \left\{ \phi_0 \in H^{1,2}(\Omega) \,\middle|\, \phi_0|_{\inflowBoundary} = 0 \right\}$ in space and time $
\mathcal{V}=\left\{\phi\,\middle|\,\phi(t) \in H^{1,2}_{\inflowBoundary}(\Omega) \,\forall\, t \in (0,T)\right\}$, cf.~\cite{Villa.2021}.
The weak formulation is then given by
\begin{equation}
\int_{0}^{T} \int_{\Omega}
\bigl(u_t + \velocity \cdot \nabla u\bigr)\, \phi
\;+\;
\kappa\, \nabla u \cdot \nabla \phi~\mathrm{d}\Omega~\mathrm{dt}
\;+\;
\int_{\Omega} u(0,\cdot)\, \phi_0~\mathrm{d}\Omega
=
\int_{0}^{T} \int_{\Omega} \parameterContinuous \, \phi~\mathrm{d}\Omega~\mathrm{dt}
\;+\;
\int_{\Omega} \parameterInitial \, \phi_0~\mathrm{d}\Omega,
\end{equation}
$\forall \phi \in \mathcal{V}, \phi_0\in H^{1,2}_{\inflowBoundary}(\Omega)$. Here, the parameter $m$ enters the formulation solely through the right-hand side.}
 Our goal is to discretize this formulation. To this end, we introduce a finite-dimensional subspace $\ansatzSpace \subset H_{\inflowBoundary}^{1,2}(\Omega)$ of continuous Lagrange nodal basis functions, defined by $\ansatzSpace = \text{span}\{\phi_1, \dots, \phi_\ndof\}$. To approximate the $L^2$-norm, the mass matrix $M \in \mathbb{R}^{\ndof \times \ndof}$ is introduced, with entries $M_{ij} := \int_{\Omega} \phi_i(\x) \phi_j(\x)~\mathrm{d}\tcr{\Omega(x)}$.  An identification between the $L^2$-functions, their norm and the discrete elements is required. For this purpose, let the map $\FEMi: \mathbb{R}^\ndof \rightarrow \ansatzSpace$ be defined by $\FEMi(a) = \sum_{i=1}^\ndof a_i \phi_i$. Together with the scalar product $\langle a, b \rangle_M := \langle M a, b \rangle_{\mathbb{R}^\ndof}$, the map $I$ becomes an isometry, i.e., 
$ \langle \FEMi(a), \FEMi(b) \rangle_{L^2(\Omega)} = \langle a, b \rangle_M$. Similarly, the stiffness matrix $K_{ij} = \int_{\Omega} \langle \nabla \phi_i(\x), \nabla \phi_j(\x) \rangle~\mathrm{d}\tcr{\Omega(x)}$ and the skew-symmetric matrix $V_{ij} = \int_{\Omega} \phi_i(\x) \langle \nabla \phi_j(\x), \velocity \rangle~\mathrm{d}\tcr{\Omega(x)}$ are defined. Additionally, the well-established SUPG-stabilization technique~\cite{Brooks.1982}, requires the matrix
$S_{ij} = \int_{\Omega} \langle \nabla \phi_i(\x), \velocity \rangle \left(\langle \nabla \phi_j(\x), \velocity \rangle -\kappa \Delta \phi_j(\x)  \right)~\mathrm{d}\tcr{\Omega(x)}$. To account for local characteristics of the initial boundary value problem, the stabilization parameter $\tau$ is defined as
$\tau = \min \left( \fracSolidus{h_E^2}{2 \kappa}, \fracSolidus{h_E}{\| \velocity \|} \right)$, where $h_E := \sup_{x, y \in E} |{x} - {y}|$ denotes the diameter of a finite element $E$. A more advanced metric-based definition of the stabilization parameter is provided in~\cite{Danwitz.2023}. 

The transient problem is addressed using an implicit Euler time-stepping algorithm~\cite[Equation (10.25)]{Elman.2014}. We use the approximation $\vec{{u}}_{t} \approx \fracSolidus{\left({\up} - {\un}\right)}{\Delta t}$ for time instances $(0, \Delta t, \dots, T = \Delta t \, n_T)$ to find a solution in discrete space-time $\bigoplus_{i=0}^{n_T} \ansatzSpace$ of the equation
\begin{align}
\label{eq:weakAD}
(M + \Delta t V + \Delta t \kappa K + \Delta t \tau S + \tau V^T){\up} = (M + \tau V^T)({\un} + \Delta t \, \vec{m}_{\text{C}}^n),
\end{align} 
with the initial condition $\vec{u}_{n=0} = \vec{m}_{\text{I}}$. The space \( \bigoplus_{i=0}^{n_T} \ansatzSpace \) represents the temporal product of the spatial ansatz space over all time steps. Function composition gives the discrete linear operator $\pts^h : D^h \rightarrow \bigoplus_{i=0}^{n_T} \ansatzSpace$ for the chosen time discretization, with the parameters projected into their corresponding finite element representations.

The dual problem for the forward equation has already been derived in \eqref{eq:adjoint_equation}. \tcr{To solve this terminal-value problem, we apply the change of variables $t = T - \bar t$, which transforms the problem \eqref{eq:adjoint_equation} into a well-posed
advection--diffusion problem with a prescribed initial condition.
The resulting problem can then be discretized using the same
methodology as for the forward equation. Under this transformation, only the direction of advection is reversed, i.e., \(-\velocity \). See, for example,
\cite{Givoli.2014} in the context of the wave equation for a detailed explanation. } To obtain the dual or adjoint problem in discrete form, two \tcr{ different} approaches can be followed. \tcr{We pursued the approach of taking the dual of the stabilized primal problem \eqref{eq:weakAD}, which leads to the following equation
\begin{align}
\label{eq:dis_adj}
(M + \Delta \bar t V^T + \Delta \bar t \kappa K + \Delta \bar t \tau S^T + \tau V) {\pn} &= (M + \tau V) {\pp} + M \vec{\misfit}^{n+1}\,,
\end{align}
with initial \( \adjoint_{n=0} = 0 \).} This variant leads to the discrete operator \( \mta^h : \mathbb{R}^{\Nobservations} \rightarrow \bigoplus_{i=0}^{n_T} \ansatzSpace \). 

In the alternative approach, the adjoint problem \eqref{eq:adjoint_equation}, which becomes a transport problem \tcr{after the transformation \(t = T - \bar t\)}, is stabilized and discretized\tcr{, and the differences in numeric accuracy will be treated in the remainder of this section}.
\begin{figure}
\centering
\includegraphics[width=0.6\linewidth]{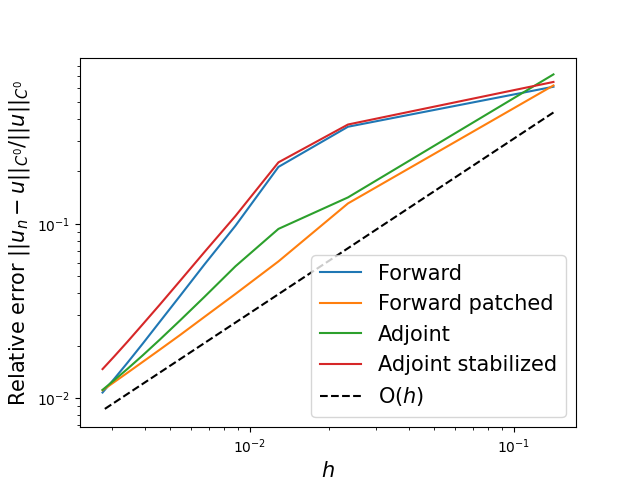} 
\caption{Asymptotic linear convergence of forward problem, patched forward problem, adjoint problem and stabilized adjoint problem of test case with Dirac initial condition $\delta_{\x_s}$ on a unit square with the analytical solution $G(\x,t)= \fracSolidus{\exp(-\fracSolidus{\norm{\x-\x_s-\velocity t}^2}{4 t \kappa})}{(4 \pi t \kappa)}$ in $[0,T] \times \R^2$ for $\x_s=[0.5,0.5],\,\kappa=\SI{0.001}{\square\metre\per\second}, \fracSolidus{h}{\Delta t} = \text{const.}$ and $\velocity=[\SI{0.1}{\metre\per\second},\SI{0.1}{\metre\per\second}]$ evaluated at $t=\SI{1}{\s}$.}
\label{fig:convergence}
\end{figure}
For the discretization, linear Lagrange finite elements of first order are used throughout the following. This choice ensures good numerical accuracy, provided that the discretization is sufficiently fine. Examining the term for \( S_{ij} \) in the SUPG-stabilization, it is evident that the diffusion term \( -\kappa \Delta \phi_i \) is not included in the residual for first-order elements. To improve accuracy, the second derivatives can be "patched." Specifically, Eq.~\eqref{eq:weakAD} is initially solved under the condition \(\Delta {\un} = 0\). Subsequently, the derivatives of the resulting solution \({\un}\) are projected onto a higher-order function space. Subsequently, the patched equation is solved using the approximate second derivatives of \( \Delta {\un} \). This approach enhances the accuracy of the forward problem, but does not improve the observed asymptotic linear convergence in the $C^0$-norm, see \autoref{fig:convergence}. In the adjoint problem (Eq.~\eqref{eq:dis_adj}), which is the dual problem of the SUPG-stabilized forward problem, the sensor misfit is not incorporated into the SUPG-stabilization, \tcr{ which would be the case if the adjoint problem \eqref{eq:adjoint_equation} were SUPG-stabilized after performing the change of variables. Despite this omission, the convergence behavior does not deteriorate compared to the "patched" fully SUPG-stabilized adjoint problem, see \autoref{fig:convergence}. For this reason, we have pursued the approach of using the dual problem of the SUPG-stabilized primal problem, which leads to Eq.~\eqref{eq:dis_adj}.}

To connect primal and dual \tcr{equation} also a discretized form of the observation operator $\obsO$ is required. Therefore, the discrete solution $\FEMi({\un}) = \sum_{i=1}^\ndof u^n_i \phi_i$ must be evaluated on the grid as observation points do not necessarily coincide with finite element nodes, i.e., $\x_i^{\text{obs}} \notin \{\x_1, \dots, \x_\ndof\}$. This process can be achieved using barycentric interpolation, as described in \cite[Problem 1.3]{Elman.2014}. For the sake of completeness, the method is briefly summarized for triangular and tetrahedral meshes. The interpolation of the function $u$, given in a discrete form, is discussed. There exists a triangle $\{\x_{l_1}, \x_{l_2}, \x_{l_3}\}$ or a tetrahedron $\{\x_{l_1}, \x_{l_2}, \x_{l_3}, \x_{l_4}\}$ for $\Omega \subset \mathbb{R}^{n_{\mathrm{sd}}=2}$ or $\subset \mathbb{R}^{n_{\mathrm{sd}}=3}$ containing $\x_i^{\text{obs}}$. The interpolation of the discrete function $\FEMi({\un})\observationforequation$ in a barycentric coordinate system is given by $\sum_{j} u_{l_j} a_j$, with $0 \leq a_{j} \leq 1$, $\sum_{j} a_j = 1$ and $\x_i^{\text{obs}} = \sum_j a_{j} \x_{l_j}$. In case $t_i^{\text{obs}}\neq k \Delta t$ for $k \in \mathbb{N}_0$ and $k \leq n_T$, interpolation between two time steps is needed \tcr{additionally} to construct the \textit{discrete observation operator} $\obsO^h : \bigoplus_{i=0}^{n_T} \ansatzSpace \rightarrow \mathbb{R}^{\Nobservations}$ and the \textit{discrete parameter-to-observable map} $\pto^h : \ansatzSpace \rightarrow \mathbb{R}^{\Nobservations}$, where $\pto^h = \obsO^h \circ \pts^h$. Moreover, using the barycentric coordinate system, the discrete Dirac distribution is defined as
\begin{equation}\label{eq:discret_dirac}
\vec{\delta}_{\x_i^{\text{obs}}} = M^{-1} \sum_{j} a_j \stdbase_{l_j} \text{, where } a_j \text{ is a barycentric representation of } x_i^{\text{obs}}.
\end{equation}
Now, the right-hand side of the adjoint problem~\eqref{eq:adjoint_equation} can now also be formulated with
$\misfit=\misfiti[i]\vec{\delta}_{\observationforequation}\in \bigoplus_{i=0}^{n_T} \ansatzSpace,$
for a given misfit vector $\misfit \in \mathbb{R}^{\Nobservations} $. \tcr{Concerning the discretization of the unknown sources and shape functions, we replace $\som(\Omega)$ by the cone of Dirac-Delta functionals supported in the grid nodes,
\begin{equation} \label{eq:discmeasures}
    \som_h=\operatorname{cone} \left\{\,\delta_{x_j}\;|\;j=1,\dots, \ndof\, \right\}= \left\{\,\sum^\ndof_{j=1} \lambda_j\delta_{x_j} \;|\; \sourceIntensity \in \R^{\ndof}_{\geq 0}\,\right\}.
\end{equation}
Subsequently, the discretized version of the radial basis function $\ansatzSources^1$ corresponds to the standard $L^2$-projection onto the respective elements. The discrete realization of the second source type, $\ansatzSources^2$, is obtained by solving $(\ellipI M + \ellipLaplace K + \beta B)\vec{m} = \vec{\delta}_{\x}$ where $B$ with $B_{ij} := \int_{\partial \Omega} \phi_i(\x) \phi_j(\x) \, d\tcr{\partial \Omega(\x)}$ corresponds to the Robin boundary condition. Finally, the third source type, $\ansatzSources^3$, does not require further discretization. Afterwards, \autoref{alg:pdap} is applied to the discretized problem where the new candidate locations are computed by
    \begin{equation*}
        \widehat{\x}_k^{\text{I}}\in \argmax_{\x \in \left\{x_1, \dots, x_\ndof\right\}} \varphi^k_{\text{I}}(\x), \quad \widehat{\x}_k^{\text{C}}\in \argmax_{\x \in \left\{x_1, \dots, x_\ndof\right\}} \varphi^k_{\text{C}}(\x)
\end{equation*}
which requires $\mathcal{O}(\ndof)$ operations.}

\tcr{In summary, this leads to two levels of approximation: On the one hand side, the finite element discretization to simulate the forward and adjoint problem. On the other hand, the discretization of the admissible set, \eqref{eq:discmeasures}. Since, in general,  optimal solutions of the limit problem \eqref{eq:sparseObjective} are not supported in grid nodes, the latter leads to clustering effects, i.e. an optimal Dirac-Delta functional in the interior of a triangle/tetrahedron is approximated by various sources in the adjacent grid nodes. In a subsequent heuristic post-processing step, we merge the latter into a single source using barycentric coordinates. This is illustrated in \autoref{fig:example_tb_baseline}.}
\tcr{
\begin{remark}
For piecewise linear and continuous finite elements as well as $\parameterInitial= \ansatzSources^3$, there holds
\begin{equation*}
    \max_{\x \in \left\{x_1, \dots, x_\ndof\right\}} \varphi_{\text{I}}(\x)=\max_{\x \in \Omega} \varphi_{\text{I}}(\x)
\end{equation*}
since $\varphi_{\text{I}}=\mta(\misfit)(0, \cdot) \in \ansatzSpace$ in this case. Hence, the admissible set is implicitly discretized due to the approximation of the PDE which is also referred to as variational discretization, see, e.g., \cite{Casas.2012}. Note, however, that this does not alleviate the mentioned clustering phenomena.
\end{remark}}

\section{Numerical Results}\label{sec:numerical}
Numerical results evaluate the performance of the presented algorithm on four test case geometries, namely on a two building benchmark geometry in 2D and 3D, on a campus geometry with wind in inflow-outflow configuration and on a 2D cut of a chemistry plant site. If not stated otherwise, numeric values of physical quantities are reported in standard SI units. 
\subsection{Benchmark geometry representing two buildings}
\begin{figure}
\begin{subfigure}{0.48\textwidth}
\includegraphics[width=0.9\linewidth]{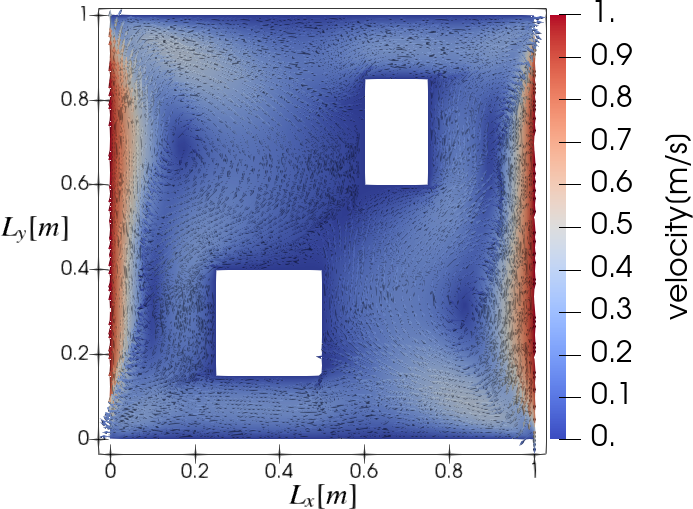}
\end{subfigure}
\begin{subfigure}{0.48\textwidth}
\includegraphics[width=0.86\linewidth]{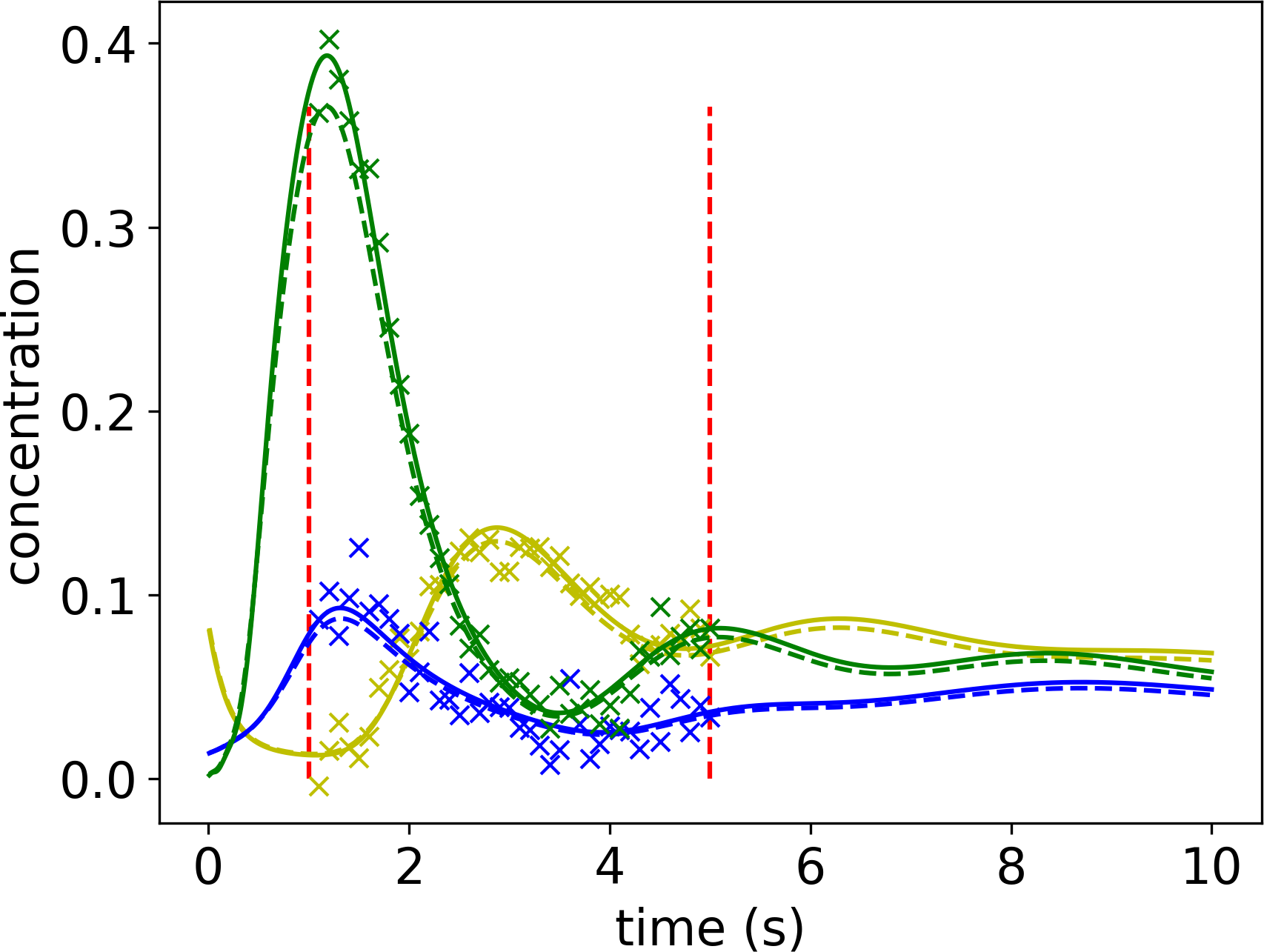}
\end{subfigure}
\caption{PDAP-algorithm input. Wind vector field (left) and noisy measurements at 3 sensors from $t=\SI{1}{\s}$ to $t=\SI{5}{\s}$ (right, marked by crosses), along with exact concentration values (full line), and PDAP-reconstruction (dashed line)}
\label{fig:example_tb_wind}
\end{figure}
First, a widely used benchmark geometry is employed to evaluate the algorithm and its implementation. The geometry consists of a unit square with two rectangular obstacles ("buildings"), and the double-lid-driven-cavity setup is utilized to create a dynamic flow scenario~\cite{Villa.2021}:
\begin{equation}\label{eq:INS_strong}
    \begin{aligned}
        -\nu \nabla^2\bm{\wind} + \bm{\wind} \cdot\nabla\bm{\wind} + \nabla \rho &= 0 \quad &\text{in} \quad &\Omega, \\
        \nabla \cdot \bm{\wind} &= 0 \quad &\text{in} \quad &\Omega, \\
        \bm{\wind} &= {\vec{\bm{g}}} \quad &\text{on} \quad &\Gamma_D. \\
        \end{aligned}
\end{equation}
In this example of the incompressible Navier-Stokes equations, $\rho$ represents the pressure, and $\nu = \fracSolidus{1}{50}$ denotes viscosity, resulting in a Reynolds number of 50. The Dirichlet boundary conditions are specified as ${\vec{g}} = {\bm{e}}_2$ on the left wall of the domain and ${\vec{g}} = -{\bm{e}}_2$ on the right wall. No-slip boundary conditions are applied on the remaining boundaries. The wind vector field for this benchmark example is depicted in \autoref{fig:example_tb_wind}, see also~\cite{Petra.2011}. To reproduce reference results~\cite{Wu.2023, Wogrin.2023}, the radial basis function $\ansatzSources^1(\cdot, 0.26, \cdot)$ is chosen as the shape function $\ansatzSourcesInitial$, and the initial condition is defined as $u_0^{\text{rbf}} = \parameterInitial(\measureInitial[{\xbr}^{\text{I}}, {\sourceIntensity}^{\text{I}}])$, with ${\xbr}^{\text{I}} = [(0.35,\, 0.7)]$ and ${\sourceIntensity}^{\text{I}} = [1.0]$ (based on Eq.~\eqref{eq:gauss_blob}). The diffusion constant $\kappa = 0.001\,\si{\square\metre\per\second}$ has been adjusted accordingly. Synthetic measurement data $\measurement$ is generated with additive white noise following $\mathcal{N}(0, \sigma^2 \operatorname{Id})$. In this case, a relative noise of \( 3\% \) was multiplied by the maximum measured concentration to determine the variance \( \sigma^2 = (0.03 \cdot \max_{i \in \measurementSpace} u\observationforequation)^2 \) for the perturbation. This corresponds to a signal-to-noise ratio of \( \text{SNR} \approx 33.3 \). The resulting measurement input to the PDAP algorithm is visualized in \autoref{fig:example_tb_wind} (right). With \( \alpha = 1000 \) as the regularization parameter (in Eq.~\eqref{eq:finitealgo}), only eight iterations are required to reconstruct the solution, which corresponds to eight forward solutions of~\eqref{eq:forward_equation}, and eight adjoint solutions of~\eqref{eq:adjoint_equation}.
\begin{figure}
\begin{subfigure}{0.32\textwidth}
\includegraphics[width=1.0\linewidth]{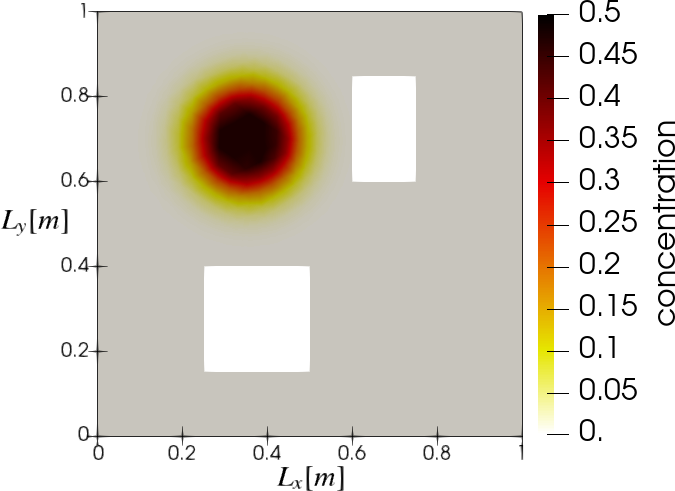}
\end{subfigure}
\begin{subfigure}{0.32\textwidth}
\includegraphics[width=0.94\linewidth]{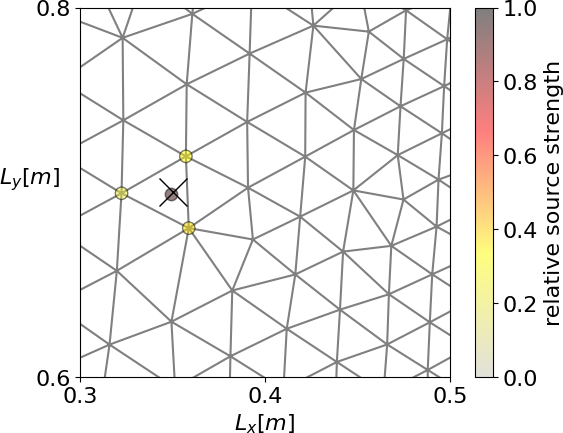}
\end{subfigure}
\begin{subfigure}{0.32\textwidth}
\includegraphics[width=1.0\linewidth]{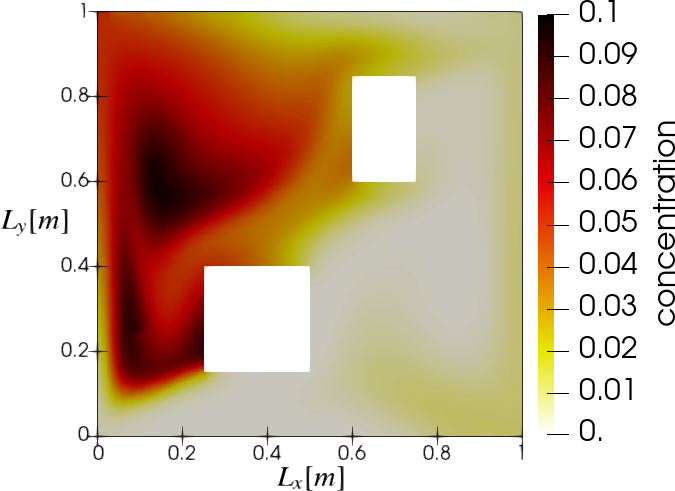}
\end{subfigure}
\caption{PDAP-algorithm output for baseline simulation. Reconstructed initial condition (left), identified parameter, namely, reconstructed raw sources $ \bar{\lambda}^{\text{I}}_i \delta_{\bar{\x}^\text{I}_i}$ (yellow dots) and post-processed source $\mu_{\text{post}}={\lambda}_{\text{post}}\delta_{{\x}_\text{post}}$ (dark dot) next to actual source center (black cross, middle), predicted concentration field at $t=\SI{5}{\second}$ (right). Ground truth in~\autoref{fig:example_tb_single}.}
\label{fig:example_tb_baseline}
\end{figure}
 The results displayed in \autoref{fig:example_tb_baseline} accurately reproduce the outcomes of established approaches that utilize $L^2$-regularization~\cite{Villa.2021,Petra.2011,Wu.2023,Alexanderian.2014, Alexanderian.2018, Attia.2018,Wogrin.2023}. Note that both presented snapshots of the computed concentration fields are in the extrapolation regime, i.e., at time instances for which no measurement data are available.
 \begin{table}[ht]
    \caption{Robustness of PDAP algorithm with respect to challenging measurements conditions in terms of  Euclidean distance from actual source center $\x_s$ to algorithmically reconstructed source location ${\x}_\text{post}$}
    \label{tab:pdap_robustness}
    \centering
    \begin{tabular}{l c c c}
    \toprule
        Measurement setting & original & reduced &$\norm{{\x}_\text{post}-\x_s}_{\Omega}$ \\
        \midrule
        Baseline simulation (\autoref{fig:example_tb_baseline}) & \multicolumn{2}{c}{ } & \SI{0.001}{\m}\\
         Reduced signal-to-noise ratio $\text{SNR}$ &  $\approx 33.3 $&$ \approx  6.7$  &  \SI{0.003}{\m}\\
         Reduced number of sensors $\numberOfSensors$ &   9 &   3  & \SI{0.013}{\m} \\
         Reduced sampling rate $f_s$ & $\SI{10}{\hertz} $ & $ \SI{2}{\hertz}$  & \SI{0.019}{\m}\\
         Reduced sampling time $T_s$ & $[\SI{1.0}{\s},\SI{5.0}{\s}] $&$ [\SI{2.0}{\s},\SI{3.0}{\s}]$  &\SI{0.017}{\m}\\
         Single moving sensor (\autoref{fig:example_baseline_moving_initial})& \multicolumn{2}{c}{ } & \SI{0.050}{\m}\\
         \bottomrule
    \end{tabular}

\end{table}
 
\begin{figure}
\begin{subfigure}{0.48\textwidth}
\includegraphics[width=0.86\linewidth]{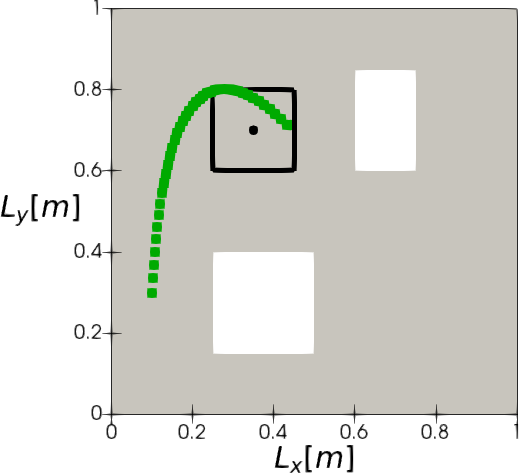}
\end{subfigure}
\begin{subfigure}{0.48\textwidth}
\includegraphics[width=1.0\linewidth]{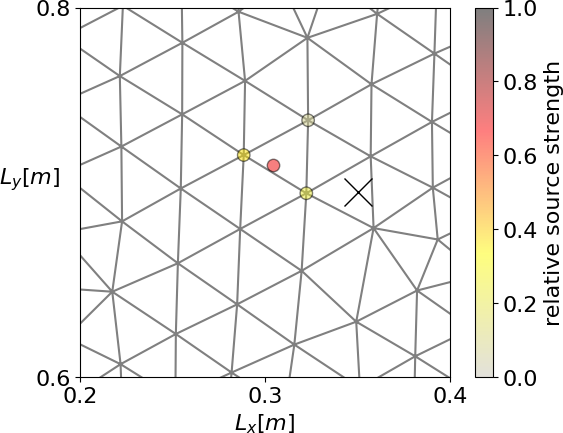}
\end{subfigure}
\caption{Single moving sensor experiment. Trajectory of the moving sensor $\gamma(t)= [0.6,0.5] + t\,[-\fracSolidus{11}{60},\fracSolidus{83}{300}] + t^2\,[\fracSolidus{1}{60},-\fracSolidus{19}{300}]$ for $t^{\text{obs}} \in \{1.0s,1.1s,...,5s\}$ (green), source location and excerpt (black, left), predicted raw sources $\lambda_i^{\text{I}}\,\delta_{\x^{\text{I}}_i}$ (yellow dots) and derived source $\mu_{\text{post}}$ (red dot) near ground truth source center (black cross, right)}
\label{fig:example_baseline_moving_initial}
\end{figure}
Additional experiments are conducted to demonstrate the robustness of \autoref{alg:pdap} under challenging measurement conditions, such as reducing the signal-to-noise ratio, the number of sensors, the sampling rate, and the sampling time. Furthermore, we emulate a moving sensor (\autoref{fig:example_baseline_moving_initial}). The results of these numerical experiments are reported in \autoref{tab:pdap_robustness}. The obtained accuracy of source identification remains acceptable under all five considered perturbations. Thus, we conclude that the algorithm is robust with respect to reduced measurement quality in this test case. Moreover, a single moving sensor with predefined trajectory is sufficient to identify the contaminant source with acceptable accuracy. 
 
\begin{figure}
\centering
\begin{subfigure}{0.49\textwidth}
\includegraphics[width=1.0\linewidth]{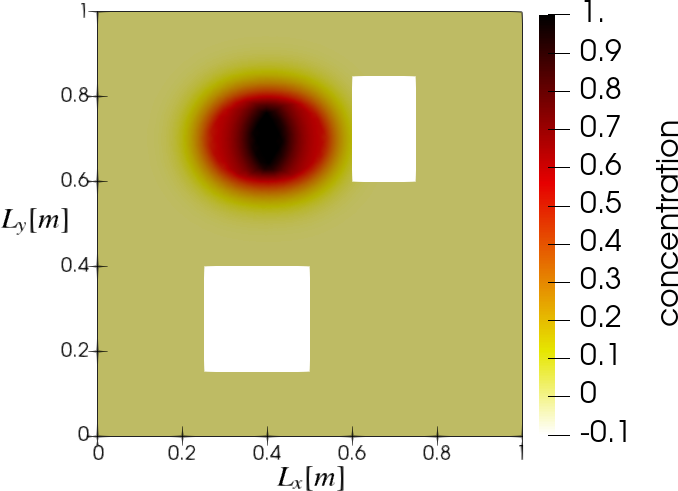}
\caption{}
\end{subfigure}
\centering
\begin{subfigure}{0.49\textwidth}
\includegraphics[width=1.0\linewidth]{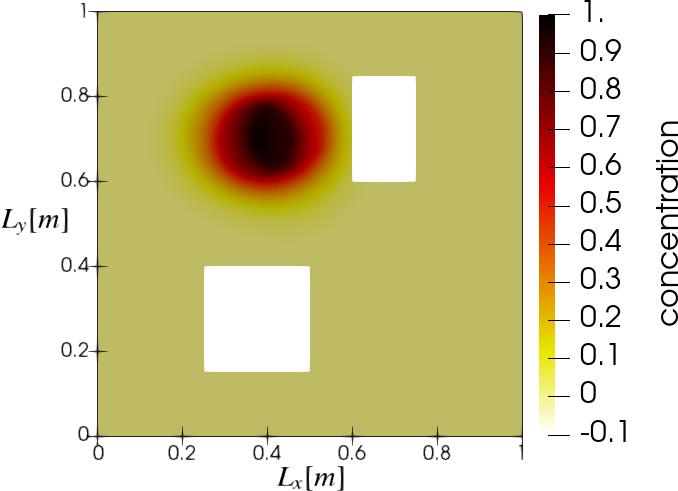}
\caption{}
\end{subfigure}
\centering
\begin{subfigure}{0.49\textwidth}
\includegraphics[width=1.0\linewidth]{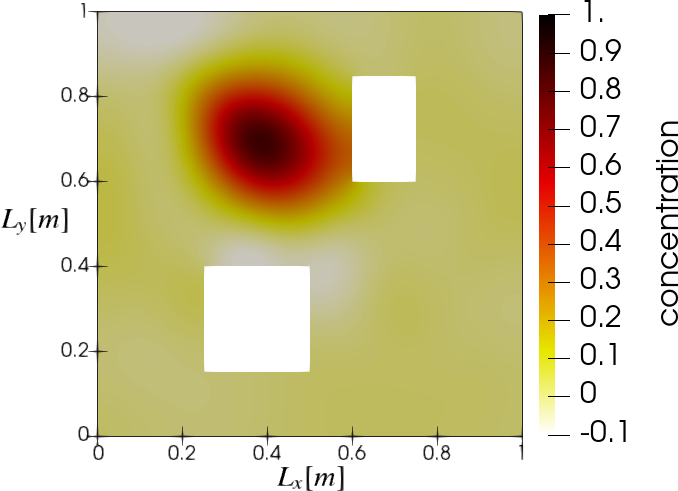}
\caption{}
\end{subfigure}
\centering
\begin{subfigure}{0.49\textwidth}
\includegraphics[width=0.93\linewidth]{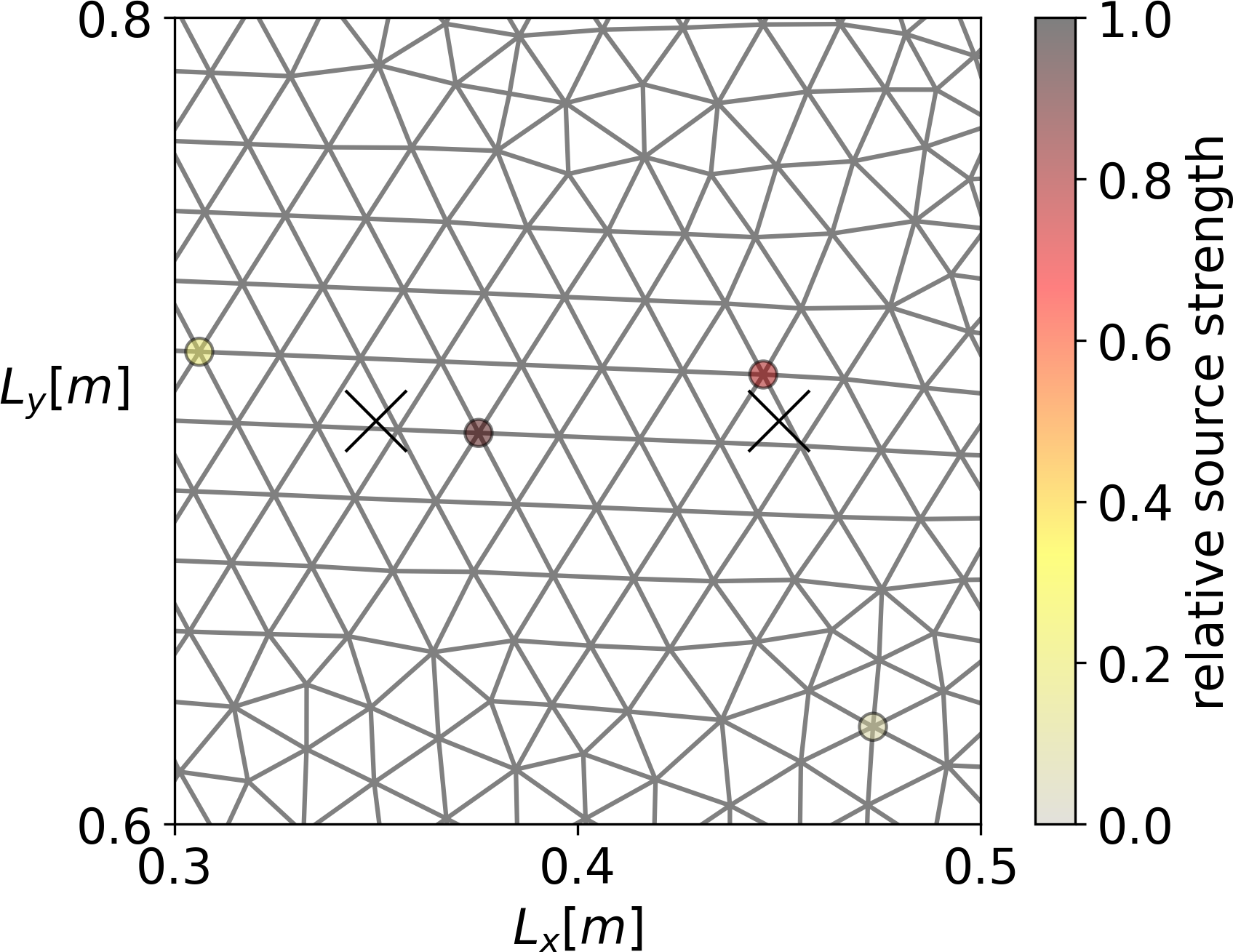}
\caption{}
\end{subfigure}
\caption{Two overlapping contaminant sources. \tcr{Ground truth initial condition $u_0^{\text{double}}$ (a), PDAP-reconstructed contaminant sources (b), $L^2$-reconstruction (c) and predicted raw sources $\sourceIntensityi[i]^{\text{I}}\delta_{\x^{\text{I}}_i}$}(d)}
\label{fig:example_tb_double}
\end{figure}
To distinguish our method from established $L^2$-approaches, we propose a more challenging initial condition with two overlapping sources $u_0^{\text{double}}= \parameterInitial(\measureInitial[{\xbr}^{\text{I}}, {\sourceIntensity}^{\text{I}}])$, with ${\xbr}^{\text{I}}=[(0.35,0.7),(0.45,0.7)]$ and ${\sourceIntensity}^{\text{I}} = [1.0,1.0]$, to be identified. Here, a higher signal-to-noise ratio $\text{SNR}=100$ was used and the regularization parameter $\alpha=500$ was adapted to the specific problem. The results are shown in \autoref{fig:example_tb_double}. It is observed that the PDAP-algorithm is able to separate the two centers, while the established $L^2$-approach recovers a single smooth function and cannot provide the insight that the contaminant stems from two separate sources. \tcr{For the \(L^2\)-reconstruction, we adhere to the provided numerical example in \cite{Villa.2021} and employ the same prior distribution, given by $\mathcal{A} = (8~I - \Delta)$.}

\begin{remark}
\tcr{For the numerical examples, we set the tolerance of the Newton algorithm to a very small value, namely $\text{tol}_{\text{Newt}} = 10^{-10}$. The termination of \autoref{alg:pdap} is controlled by  $\text{tol}$, and an error estimate relative to $J(0,0)=\fracSolidus{1}{(2\,\sigma^2)}\,\| \measurement\|^2_{\R^{\Nobservations}}$ is given in \eqref{eq:residual}. In our examples, we choose $\fracSolidus{\text{tol}}{\alpha}$ between $0.01$ and $0.1$.}
\end{remark}
\subsection{Source identification in three-dimensional flow field}
The following scenario extends the two-building benchmark to a 3D geometry (\autoref{fig:3D_Sim_Single}). As initial condition, again the shape function $\ansatzSources^1(\cdot,0.25,\cdot)$ is chosen and the source is located at ${\xbr}^{\text{I}}=[(0.35, 0.7, 0.26)]$, with intensity ${\sourceIntensity}^{\text{I}} = [1.0]$. The observation pattern corresponds to a sensor array mounted with antennas on the roofs of the two cuboid buildings. The PDAD algorithm requires only 22 iterations (22 forward and 22 adjoint solutions) to accurately reconstruct the initial condition. The center of the reconstructed source is in close vicinity to the true center. A visual comparison of the original simulation and the measurement-based reconstruction is shown in \autoref{fig:3D_Sim_Single}.

\subsection{Source identification on complex campus domain with wind in inflow-outflow configuration}
\begin{figure}
\begin{subfigure}{0.48\textwidth}
\includegraphics[width=1.0\linewidth]{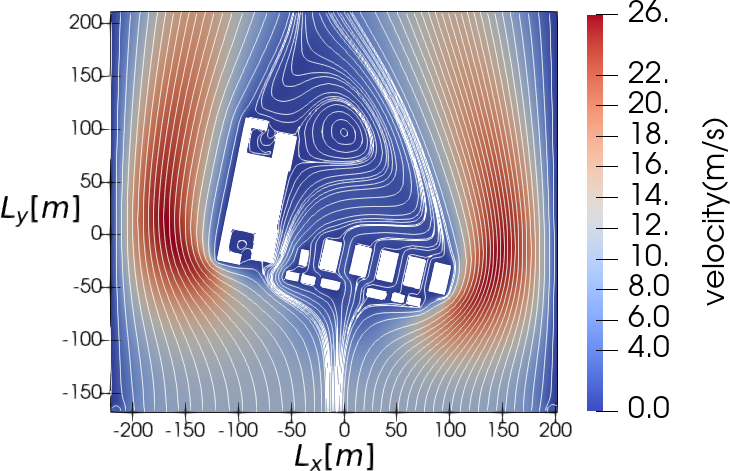}
\end{subfigure}
\begin{subfigure}{0.48\textwidth}
\includegraphics[width=0.86\linewidth]{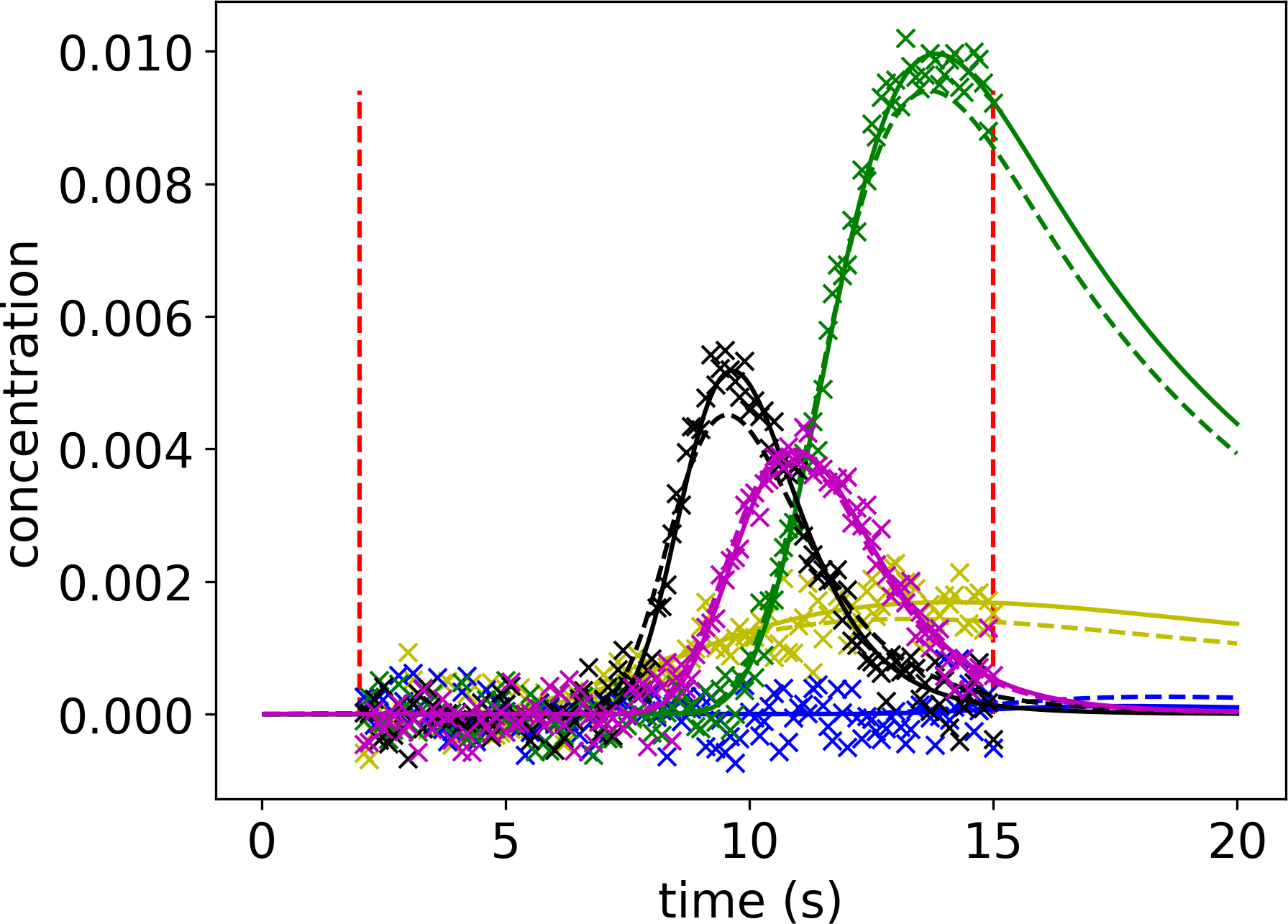}
\end{subfigure}
\caption{Campus domain input. Wind vector field (left) and noisy measurements $\measurement$ at 5 sensors (right, marked by crosses, sensor positions shown in \autoref{fig:example_campus_grid_collect}), along with exact and PDAP-reconstructed concentration (full, dashed line).}
\label{fig:example_campus_wind}
\end{figure}

To explore a more realistic scenario, a portion of the campus of the University of the Bundeswehr Munich, Germany, is selected as computational domain. The geometry comprises several building imprints directly imported from Open Street Map (OSM). Moreover, a boundary-refined finite element mesh is generated in an automated workflow~\cite{Bonari.2024}. At the southern side of the geometry, wind enters the domain resulting in the vector field shown in \autoref{fig:example_campus_wind}. For this test case, a diffusion coefficient of $\kappa = \SI{10.0}{\square \metre\per\second}$ is chosen in the transport problem \eqref{eq:forward_equation}. The observation time frame $T = [\SI{2}{\s}, \SI{15}{\s}]$ is sampled at a rate of $\SI{10}{\Hz}$, resulting in the input signal displayed in \autoref{fig:example_campus_wind}.

\subsubsection{Closed form initial condition -- comparison with \texorpdfstring{$L^2$}-approach}

\begin{table}
    \caption{Computational cost of established $L^2$-approach and PDAP-algorithm reported in terms of required PDE solutions}
    \label{tab:comp_pdap}
    \centering
\begin{tabular}{ccccccc}
\toprule
 &              &           &           & & \multicolumn{2}{c}{PDE solutions}  \\
 case & source type & \# sensors & algorithm &  parameters &  online & offline \\
\midrule
\protect\hypertarget{label_a}{a} & $u^{\text{rbf}}_0$ & 213 & PDAP & $\alpha=1000$ & 17 & -  \\
\protect\hypertarget{label_b}{b} & $u^{\text{rbf}}_0$ & 213 &$L^2$ & $\ellipI=160.0,\,\ellipLaplace=20.0$ & 54 & 1600  \\
\midrule
\protect\hypertarget{label_c}{c}  & $u^{\text{rbf}}_0$ & 5 &PDAP  & $\alpha=1500$ & 20 & -  \\
\protect\hypertarget{label_d}{d} & $u^{\text{rbf}}_0$ & 5 &$L^2$ & $\ellipI=160.0,\,\ellipLaplace=20.0$ & 16 & 120  \\
\midrule
\protect\hypertarget{label_e}{e} & $u^{\text{ellip}}_0$ & 213 & PDAP & $\alpha=1000$ & 25 & - \\
\protect\hypertarget{label_f}{f} & $u^{\text{ellip}}_0$ & 5 & PDAP& $ \alpha=1500$ & 22 & - \\
\midrule
\protect\hypertarget{label_g}{g} & $u^{\text{dirac}}_0$ & 5 & PDAP& $\alpha=2000$ & 19 & - \\
\bottomrule
\end{tabular}
\end{table}

\captionsetup[subfigure]{labelformat=empty}

\begin{figure}
\begin{subfigure}{0.32\textwidth}
\includegraphics[width=1.0\linewidth]{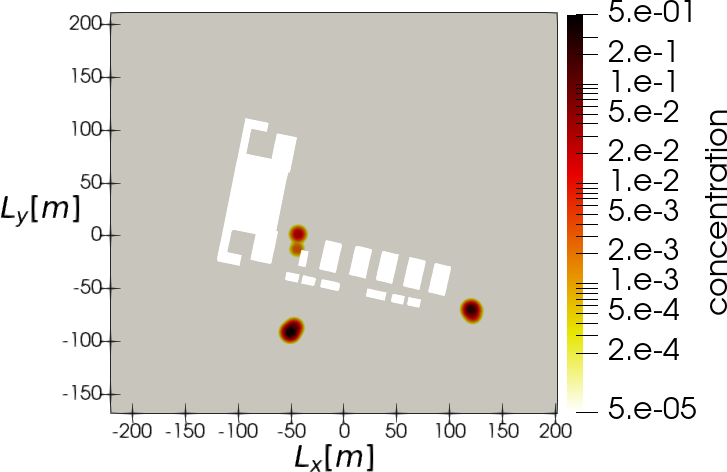}
\caption{PDAP-reconstruction, case \protect\hyperlink{label_a}{a}}
\end{subfigure}
\begin{subfigure}{0.32\textwidth}
\includegraphics[width=1.0\linewidth]{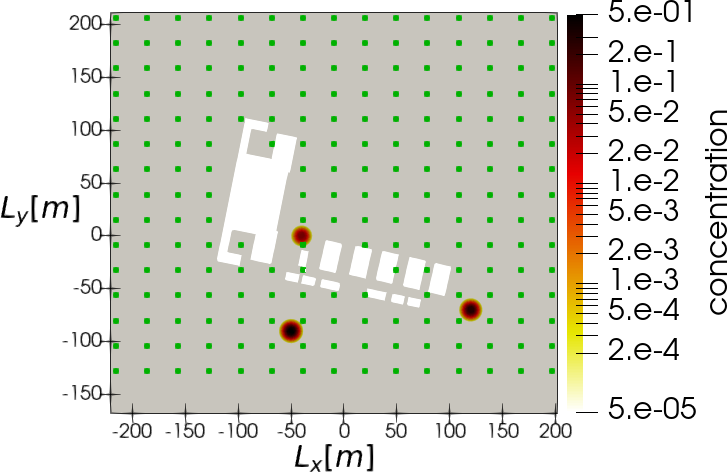}\caption{213 sensors, ground truth $u^{\text{rbf}}_0$}
\end{subfigure}
\begin{subfigure}{0.32\textwidth}
\includegraphics[width=1.0\linewidth]{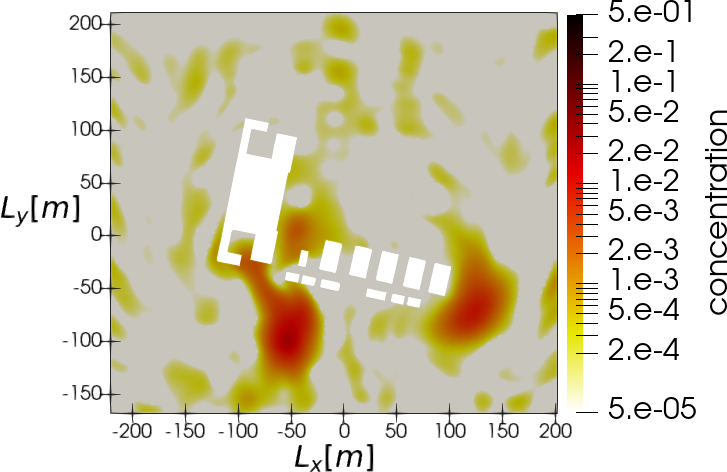}
\caption{$L^2$-reconstruction, case \protect\hyperlink{label_b}{b}}
\end{subfigure}\\
\vspace{0.5cm}
\begin{subfigure}{0.32\textwidth}
\includegraphics[width=1.0\linewidth]{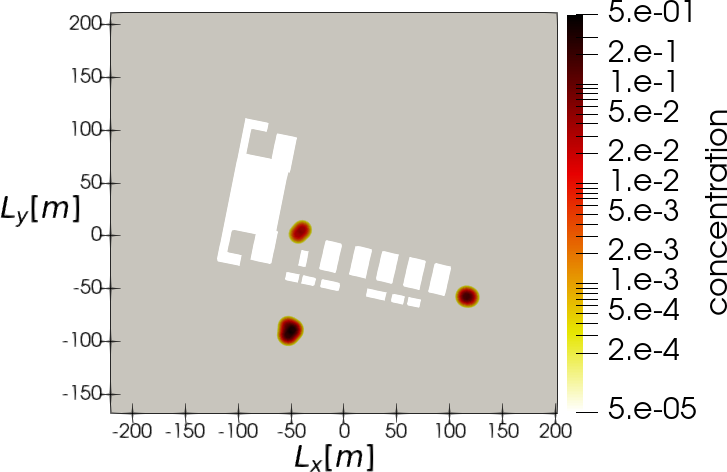}
\caption{PDAP-reconstruction, case \protect\hyperlink{label_c}{c}}
\end{subfigure}
\begin{subfigure}{0.32\textwidth}
\includegraphics[width=1.0\linewidth]{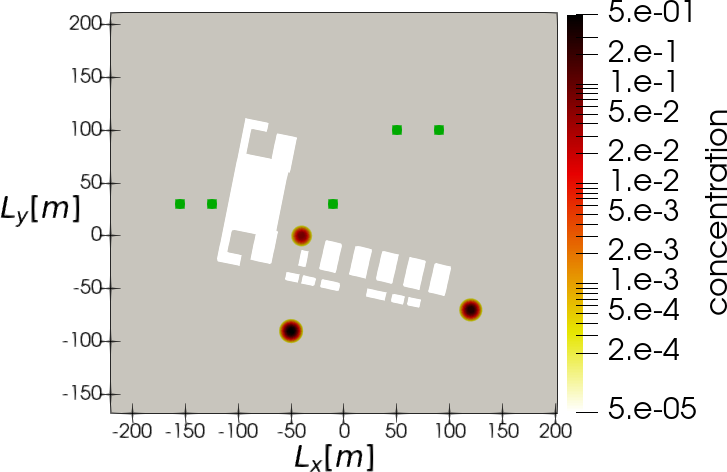}
\caption{5 sensors, ground truth $u^{\text{rbf}}_0$}
\end{subfigure}
\begin{subfigure}{0.32\textwidth}
\includegraphics[width=1.0\linewidth]{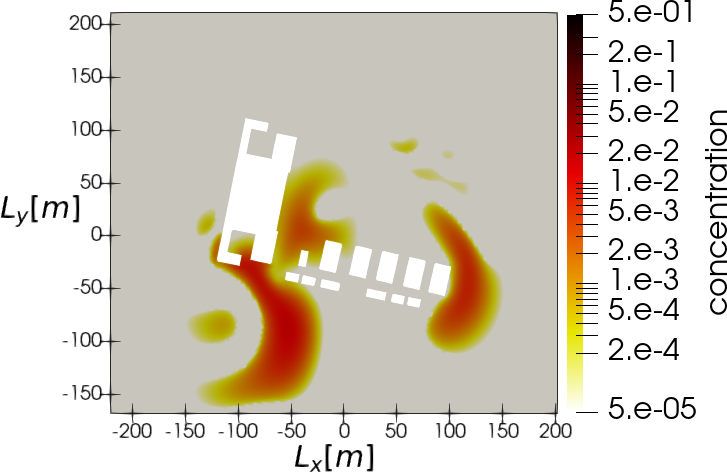}
\caption{$L^2$-reconstruction, case \protect\hyperlink{label_d}{d}}
\end{subfigure}
\caption{Comparison of PDAP-reconstruction (left) and $L^2$-reconstruction (right) of initial condition $u^{\text{rbf}}_0$ (ground truth, middle) based on measurements at sensor positions (green squares, middle) for cases listed in \autoref{tab:comp_pdap}.}
\label{fig:example_campus_grid_collect}
\end{figure}

This section elaborates how \autoref{alg:pdap} compares to an established approach with $L^2$-regularization. To achieve this, the previous initial condition is adapted to a more complex scenario in the current environment: ${u^{\text{rbf}}_0 = \parameterInitial(\measureInitial[{\xbr}^{\text{I}}, {\sourceIntensity}^{\text{I}}])}$, with ${\xbr}^{\text{I}}=[(-50, -90),(120, -70),(-40, 0)]$ and ${\sourceIntensity}^{\text{I}} = [1.0,0.5,0.1]$ and radius $r = \SI{10}{\m}$. Further algorithmic parameters are reported in \autoref{tab:comp_pdap}, along with the results for case \protect\hyperlink{label_a}{a} and case \protect\hyperlink{label_b}{b} for a relatively dense grid of 213 sensors, as well as results for a sparse sensor layout of only 5 sensors (case \protect\hyperlink{label_c}{c} and case \protect\hyperlink{label_d}{d}). Considering the combined number of PDE solutions in the online and offline phase, the presented PDAP-algorithm clearly outperforms the $L^2$-approach. Considering only online solves, an acceptable number of iterations is required in both algorithms. Notably, the number of PDAP iterations decreases for added information (more sensors) in contrast to the $L^2$-approach. \autoref{fig:example_campus_grid_collect} displays the reconstructed initial conditions and demonstrates the better accuracy of the reconstructions obtained with the PDAP-algorithm.

\subsubsection{Implicit source definition with elliptic PDE}

\begin{figure}
\begin{subfigure}{0.32\textwidth}
\includegraphics[width=1.0\linewidth]{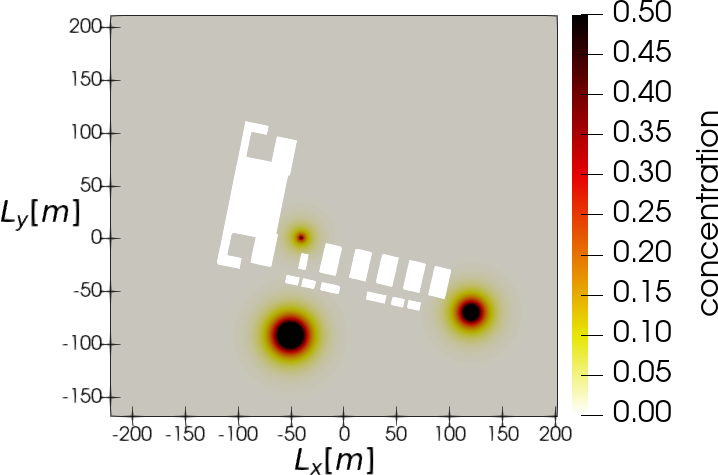}
\end{subfigure}
\begin{subfigure}{0.32\textwidth}
\includegraphics[width=0.98\linewidth]{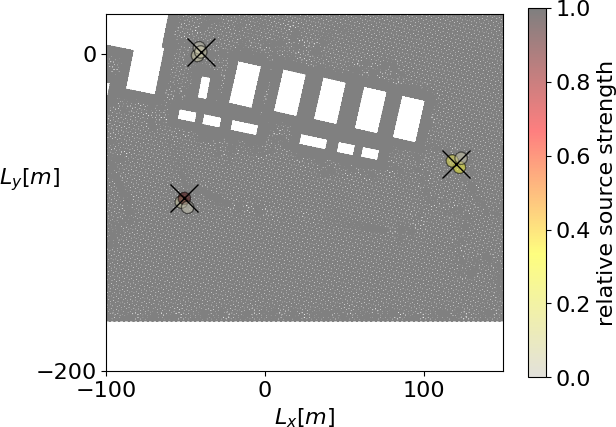}
\end{subfigure}
\begin{subfigure}{0.32\textwidth}
\includegraphics[width=0.98\linewidth]{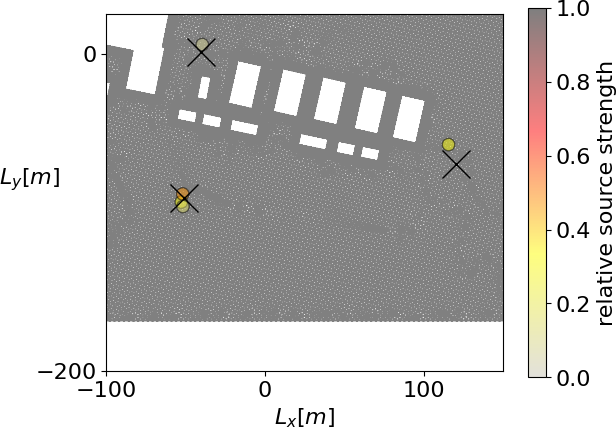}
\end{subfigure}
\caption{PDAP reconstruction of elliptic initial condition  $u^{\text{ellip}}_0$ (ground truth, left) for two sensor configurations. Reconstructed point sources based on all 213 sensors (case \protect\hyperlink{label_e}{e} in \autoref{tab:comp_pdap}, middle), and reconstruction using reduced sensor configuration (case \protect\hyperlink{label_f}{f} in \autoref{tab:comp_pdap}, right).}
\label{fig:campus_source_ellip}
\end{figure}

If the initial condition is selected as in the previous case, a convolution of the adjoint solution with the shape functions must be performed at each point in $\Omega$ to compute the dual. This process incurs a significant computational cost. To mitigate this, an alternative modeling approach for the initial condition, $\ansatzSources^2$, is considered. In this approach, the computation of $\varphi_{\text{I}}(\x)$ or $\varphi_{\text{C}}(\x)$ requires solving only a single solve of an elliptic problem (Eq.~\eqref{eq:elliptic_initial}). To motivate this, we utilize the symmetry of the operator $\mathcal{A}$ and observe that 
\begin{equation}
\int_\Omega \ansatzSourcesInitial^2(\x,z) \, \mta(\misfit)(t,z)~\mathrm{d}{\Omega(z)}=\int_\Omega \delta_{\x}(z)\,\mathcal{A}^{-1}(\mta(\misfit)(t,\cdot)) (z)~\mathrm{d}{\Omega(z)}=\mathcal{A}^{-1}(\mta(\misfit)(t,\cdot))(\x)
\end{equation}
From this identity, we deduce that 
\begin{equation}
\varphi_{\text{I}}(\x) = \mathcal{A}^{-1}(\mta(\misfit)(0,\cdot))(\x),
\qquad \varphi_{\text{C}}(\x)=\int_0^T\,\mathcal{A}^{-1}(\mta(\misfit)(t,\cdot))(\x)~\mathrm{dt}
\end{equation}
and therefore replace the convolution of the shape functions with the elliptic equation solve in Step 2 in \autoref{alg:pdap}.

In the corresponding experiment (\autoref{fig:campus_source_ellip}), the initial condition is expressed as $u^{\text{ellip}}_0 = \parameterInitial(\measureInitial[{\xbr}^{\text{I}}, {\sourceIntensity}^{\text{I}}])$, with $\ansatzSources^2$ as shape function, ${\xbr}^{\text{I}}=[(-50, -90),(120, -70),(-40, 0)]$, ${\sourceIntensity}^{\text{I}} = [1000,500,100]$, $\ellipI = 1$ and $\ellipLaplace = 100$. Parameters are selected to closely match the appearance of the initial condition $u^{\text{rbf}}_0$ above. Further details of the configuration are provided in case \protect\hyperlink{label_f}{f} and very satisfactory results are presented in \autoref{fig:campus_source_ellip}. In comparison to case \protect\hyperlink{label_a}{a} and \protect\hyperlink{label_c}{c}, a slight increase of PDAP-iterations is observed. However, this approach is still preferred as Step 2 in \autoref{alg:pdap} is greatly simplified, i.e., no convolution is required.

\subsubsection{Instantaneous point sources modeled as Dirac distributions}

\begin{figure}
\begin{subfigure}{0.48\textwidth}
\includegraphics[width=1.0\linewidth]{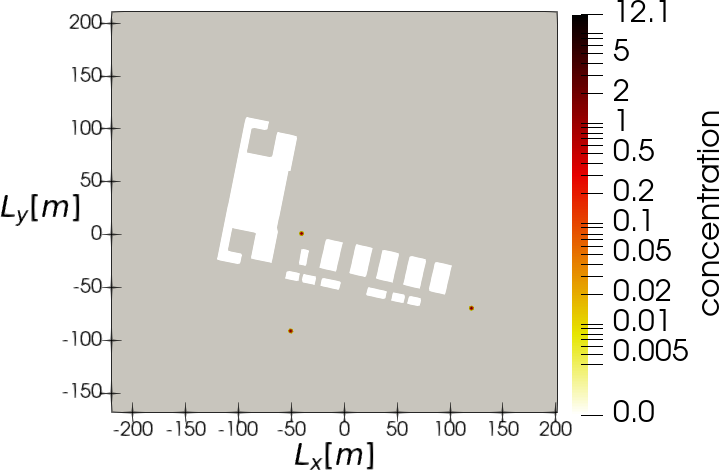}
\end{subfigure}
\begin{subfigure}{0.48\textwidth}
\includegraphics[width=.95\linewidth]{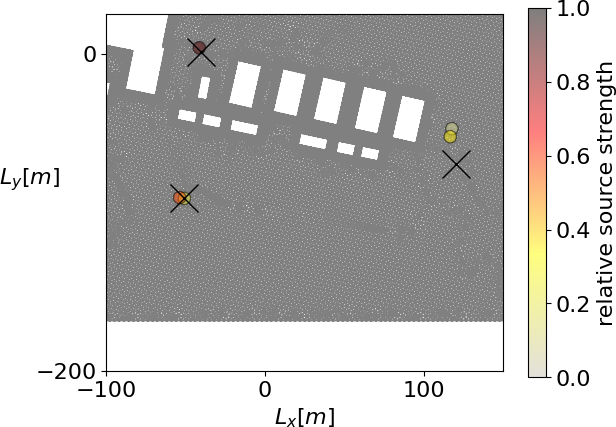}
\end{subfigure}
\caption{Campus geometry with very localized initial condition of Dirac type $u^{\text{dirac}}_0$ (left), reconstructed point sources for case \protect\hyperlink{label_g}{g} listed in \autoref{tab:comp_pdap}.}
\label{fig:example_campus_dirac_fwd}
\end{figure}

As a third example on the campus domain (\autoref{fig:example_campus_dirac_fwd}), we consider a scenario in which the initial condition \tcr{ is highly localized, taking the form of a linear combination of  Dirac-Delta functionals.} An important advantage of this formulation is that the associated dual, as defined in \autoref{thm:optimality}, coincides with the solution of~\eqref{eq:adjoint_equation}. Consequently, the convolution with the initial condition is completely omitted, which significantly simplifies the implementation and increases the computational efficiency. The use of a Dirac distribution $\ansatzSources^3 = \delta_{x_s}$ as the initial condition is also physically justified, if the area of the initially affected region is very small relative to the overall domain. In the context of a contaminant release at a chemical plant site, this can be a reasonable assumption. In the numerical experiment, the initial condition is defined as 
$u^{\text{dirac}}_0 = \parameterInitial(\measureInitial[{\xbr}^{\text{I}}, {\sourceIntensity}^{\text{I}}])$, with $\ansatzSources^3$ as shape function, ${\xbr}^{\text{I}}=[(-50, -90),(120, -70),(-40, 0)]$ and ${\sourceIntensity}^{\text{I}} = [100,100,100]$. The algorithm described in \autoref{alg:pdap} successfully identifies the sources with the dense sensor arrangement\tcr{,} and even with the reduced number of five sensors, the method maintains a satisfactory level of accuracy; see \autoref{fig:example_campus_dirac_fwd} for a visualization of the result.

\subsubsection{Continuous contaminant sources}

\begin{figure}
\begin{subfigure}{0.32\textwidth}
\includegraphics[width=1.0\linewidth]{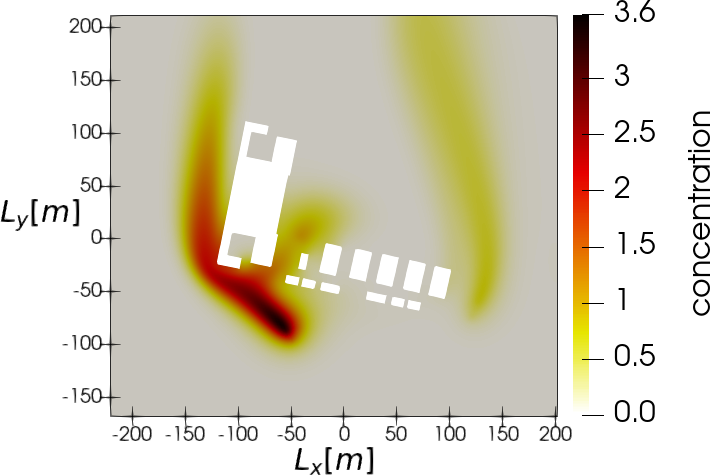}
\end{subfigure}
\begin{subfigure}{0.32\textwidth}
\includegraphics[width=1.0\linewidth]{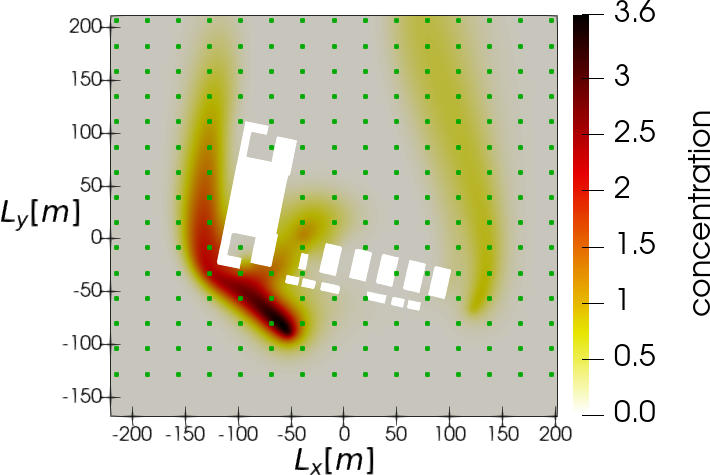}
\end{subfigure}
\begin{subfigure}{0.32\textwidth}
\includegraphics[width=1.0\linewidth]{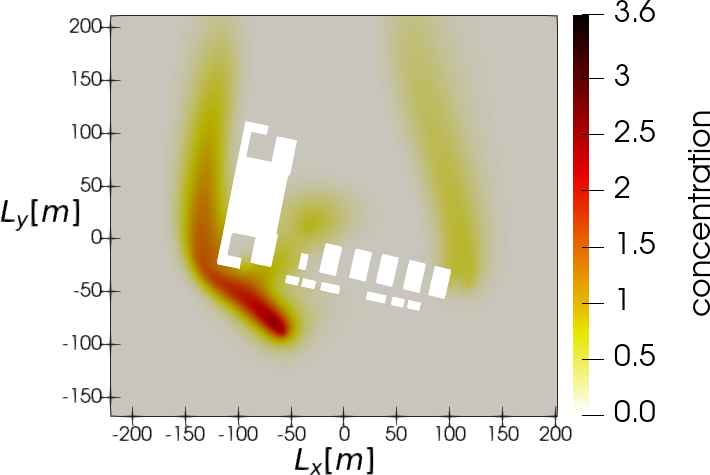}
\end{subfigure}
\caption{Identification of continuous contaminant sources. Simulated concentration field with source term $u_{\text{cont}}$ evaluated at $t=\SI{20}{\s}$ (ground truth, middle), PDAP reconstruction with 213 sensors (left), and with 5 sensors (right)}
\label{fig:example_campus_cont_fwd}
\end{figure}

The final experiment on the campus domain (\autoref{fig:example_campus_cont_fwd}) evaluates the capability of \autoref{alg:pdap} to identify continuous sources, defined as $u_{\text{cont}} = \parameterContinuous(\measureContinuous[{\xbr}^{\text{C}}, {\sourceIntensity}^{\text{C}}])$, with $\ansatzSources^2$ as shape function, ${\xbr}^{\text{I}}=[(-50, -90),(120, -70),(-40, 0)]$, ${\sourceIntensity}^{\text{I}} = [1000,500,100]$, $\ellipI = 1$ and $\ellipLaplace = 100$. As indicated by the formulation in \eqref{eq:adjoint_equation}, the adjoint solution $\adjoint=\mta(\misfit)$ may not be smooth at every time instance. However, sufficient regularity of $\varphi_{\text{C}}$ can still be ensured by selecting an appropriate observation operator $\obsO$ or applying convolution with a suitably regularized kernel. For instance $\ansatzSources^2$ provides sufficient regularity at $\varphi_{\text{C}}$ and the problem remains well-posed and solvable. In case of a densely distributed sensor network, the concentration field is again accurately reconstructed; see \autoref{fig:example_campus_cont_fwd}. However, with a sparse sensor configuration, the problem becomes more difficult. In such cases, careful sensor placement is essential to enable accurate source identification.

\subsection{Empirical robustness study of PDAP algorithm for chemistry plant site}
\begin{figure}
\begin{subfigure}{0.45\textwidth}
    \includegraphics[width=.9\linewidth]{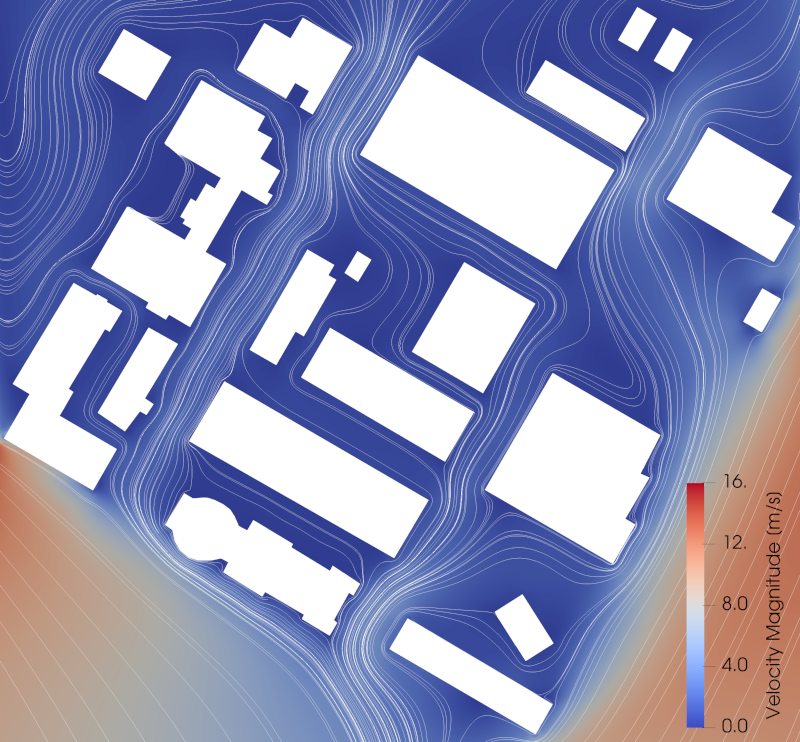}
\end{subfigure} 
\begin{subfigure}{0.53\textwidth}
    \includegraphics[width=1.00\linewidth]{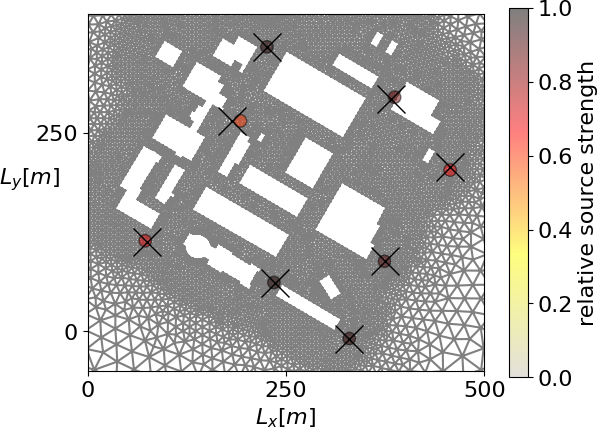}
\end{subfigure} 
    \caption{Estimated wind vector field for chemistry plant site (left), post-processed PDAP-reconstruction of sources (dots) properly indicates locations of actual sources (crosses) on  computational domain (right). Geometry adopted from~\cite{Schneider.2025}.}
    \label{fig:henkel_wind_and_source_post}
\end{figure}
The geometry of this test case is taken from the recent work of Schneider et al.~\cite{Schneider.2025}. The objective of their study is to investigate a potential application of the Emergency Response Inference Mapping (ERIMap) in collaboration with the fire brigade of a chemistry plant operator. Specifically, an accident of a tank wagon transporting chlorine is defined as test scenario. The dispersion of the released contaminant is modeled in their work using the ALOHA (Areal Locations of Hazardous Atmospheres) software, which employs a Gaussian plume model for dispersion calculations~\cite{Brusca.2016}. The accident location is assumed to be known and the primary goal of the ERIMap approach is to identify areas where individuals may have been exposed to hazardous contaminant concentrations, based on discrete sensor measurements and their associated uncertainties. This contribution extends the original approach in two significant ways. First, solving the advection–diffusion equation can provide a more realistic dispersion prediction, accounting for the influence of obstacles, such as buildings, on the contaminant transport. Second, the applicability is extended to cases with unknown source location and intensity, enabling the use of existing sensor networks for source identification. 

To estimate the wind vector field (Eq.~\eqref{eq:INS_strong}), a free-slip boundary condition was imposed on the building imprints by means of a Nitsche method as described in~\cite{Gjerde.2022}. This modeling approach results in slightly increased wind speeds in the inter-building regions, thereby providing a challenging scenario for source localization as visualized in \autoref{fig:henkel_wind_and_source_post}. The scenario is particularly difficult due to the complex interplay between diffusion and advection in the confined spaces between buildings, where flow dynamics are very heterogeneous. 
Regarding the source identification problem, eight instantaneous concentration sources are distributed across the domain. Concentration measurements were recorded at 154 sensor locations from $t_1 = \SI{2}{\s}$ to $t_2 = \SI{120}{\s}$, using a sampling rate of \SI{2}{\hertz}. As in previous experiments, synthetic measurements were generated assuming a signal-to-noise ratio of $\text{SNR} \approx 33.3$. 

In this scenario, the PDAP-algorithm required 27 iterations to achieve convergence for the results presented in \autoref{fig:henkel_wind_and_source_post}. \tcr{Based on empirical evaluation, we set the regularization parameter to $\alpha=1000$, as it provides good adaptation to the specific problem}. The computations were performed on a workstation equipped with an Intel i9-10980XE processor and took less than 2 minutes. The majority of the computation time was spent in a total of 54 PDE solves (27 forward and 27 adjoint) each with approximately 10,000 degrees of freedom. Despite the scenario's complexity, the proposed \autoref{alg:pdap} successfully identifies all eight sources with only a small number of candidate locations. Owing to the sparsity promoting regularization, these candidate locations are highly concentrated around the true source positions. For the employed sensor grid, the maximum absolute error in the prediction of the concentration field with the reconstructed sources, $u(\SI{200}{\s},\cdot)$ is bounded below $6.4\times 10^{-3}$. Due to the strong clustering, the raw results can be improved by post-processing as \tcr{in the previous section}. Candidates within a circle of \SI{40}{\meter} are projected in a post processing step to a single source.
\begin{table}
\centering
\begin{tabular}[b]{c c c c c c}
\toprule
 \multicolumn{2}{c}{Source location $\x_s$} & \multicolumn{2}{c}{Reconstruction $\x_{\text{post}}$}& $\norm{\x_{\text{post}}-\x_s}_{\Omega}$ & Estimated intensity $\sourceIntensity_{\text{post}}$\\
 \multicolumn{1}{c}{x} & \multicolumn{1}{c}{y} & \multicolumn{1}{c}{x} & \multicolumn{1}{c}{y} & \multicolumn{1}{c}{[\si{\m}]} & (truth 1.0)\\
\midrule
$329.5$&$-9.1$ &  $329.7$&$-8.9$ & 0.29 & 0.9 \\
$236.1$&$61.6$ &  $234.8$&$61.7$ & 1.33 & 1.0 \\
$382.7$&$293.4$ &  $386.9$&$295.1$ & 4.52 & 0.9 \\
$181.9$&$264.8$ &  $192.8$&$265.3$ & 10.95 & 0.6 \\
$456.6$&$207.0$ &  $456.9$&$203.4$ & 3.67  & 0.7 \\
$75.1$&$112.9$ &  $72.7$&$114.2$ & 2.79  & 0.7 \\
$226.0$&$358.3$ &  $226.0$&$358.6$ & 0.29 & 1.0 \\
$374.2$&$88.5$ &  $374.2$&$88.5$ & 0.00  & 0.9 \\
\bottomrule
\end{tabular}
\caption{Quantitative results of source localization at chemistry plant site, accompanying  \autoref{fig:henkel_wind_and_source_post}}
\label{tab:henkel_quant_source_post}
\end{table}
The quantitative results in \autoref{tab:henkel_quant_source_post} illustrate that the maximum distance from a reconstructed source to an actual source is less then \SI{11}{\meter} on the considered computational domain of about $\SI{500}{\meter}\times\SI{500}{\meter}$ and on average much smaller. Moreover, is is observed, that the reconstruction is less accurate in areas where diffusion ($\kappa= \SI{2}{\square \metre\per\second}$) dominates. All in all, \autoref{fig:henkel_wind_and_source_post} confirms that the final reconstruction provides a reliable and accurate indication for the source locations and intensities.

\section{Conclusion and Outlook}\label{sec:conclusion}
In this paper, we considered a variational regularization ansatz for the identification of an unknown number of sources from scarce measurements in contaminant dispersion applications. \tcr{The presented approach combines the theoretical insights of lifting of the problem from an unknown number of location-intensity pairs to a Radon measure on the spatial domain (\autoref{thm:optimality}) with the practical adaptation of the efficient Primal-Dual-Active-Point method (\autoref{alg:pdap}). Numerical results showcased that the presented method provides a reliable reconstruction of contaminant sources in complex environments and that the presented method outperforms state-of-the-art techniques with $L^2$-regularization.} This is particularly evident in \autoref{fig:example_campus_grid_collect}, as well as in \autoref{tab:comp_pdap}, which demonstrate that accurate source identification is possible with sparse information of few sensors. Furthermore, the presented method has the advantage that adding sensor information does not increase the computational effort in terms of PDE solves, but rather reduces it. As a result, the method scales to scenarios in which a large area is monitored with several sensors (see~\autoref{fig:henkel_wind_and_source_post}).

For future work, the requirement of the known time horizon could be lifted, i.e., the time at which the contaminant was released should be subject to identification as well. Moreover, from a computational perspective, a variety of speed-ups are possible in order to further enhance real-time capabilities of the presented approach. Many of the required PDE solves can be performed in parallel or even be precomputed. For example, by linearity, the solution of the adjoint equation \eqref{eq:adjoint_equation} in every iteration can be reduced to a single matrix-vector product if we compute and store the fundamental solutions associated to the observation points in an offline library. To further reduce computational (online) cost, suitable reduced-order models or surrogates for the parameter-to-observable map $\pto$ can be constructed using, e.g., singular value decomposition~\cite{Halko.2011}. Taking this further, integration with parameter-dependent reduced models for wind fields~\cite{Bonari.2024} could enable real-time forecasting in emergency scenarios, offering first responders fast estimates of the release locations and contaminant dispersion.

Finally, the source identification approach proposed here can be combined easily with uncertainty quantification and the optimal design of experiments~\cite{Huynh.2024}. In the context of crisis management, an assessment of the overall involved uncertainty is very helpful, e.g., to steer further measurements and support informed decision making in a Bayesian framework~\cite{Danwitz.2024, Schneider.2025, Mattuschka.03.07.2025}. Given the accurate and robust source identification, an obvious next step is to transfer the method to other applications in critical infrastructure protection (CIP), e.g., to enhance acoustic source identification with targeted applications in structural health monitoring (SHM)~\cite{Sun.2022, Scholl.2023}, or to help identifying seismic sources using an elasto-acoustic modeling and simulation approach~\cite{Mazzieri.2022}. A further methodological development might also consider sparsity promoting regularization to support contaminant source identification in water distribution networks~\cite{Jerez.2021}.

\section*{Acknowledgements}\label{sec:acknowledgements}

We would like to thank Max Winkler for discussions leading to the proper problem formulation and Johannes Pfefferer for his helpful comments on optimization strategy and implementation.

\bibliographystyle{elsarticle-num}

\bibliography{references}

@article{Flinth25,
 author = {Flinth, Axel and de Gournay, Fr{\'e}d{\'e}ric and Weiss, Pierre},
 title = {Grid is good. {Adaptive} refinement algorithms for off-the-grid total variation minimization},
 fjournal = {OJMO. Open Journal of Mathematical Optimization},
 journal = {OJMO, Open J. Math. Optim.},
 issn = {2777-5860},
 volume = {6},
 pages = {1--27},
 year = {2025},
 language = {English},
 doi = {10.5802/ojmo.39},
 zbMATH = {8032173}
}

@misc{wachsmuthwalter,
      title={No-gap second-order conditions for minimization problems in spaces of measures}, 
      author={Gerd Wachsmuth and Daniel Walter},
      year={2024},
      primaryClass={math.OC},
      doi={10.48550/arXiv.2403.12001},
}

@article{Alexanderian.2014,
 author = {Alexanderian, Alen and Petra, Noemi and Stadler, Georg and Ghattas, Omar},
 year = {2014},
 title = {A-Optimal Design of Experiments for Infinite-Dimensional Bayesian Linear Inverse Problems with Regularized {$\ell_0$}-Sparsification},
 pages = {A2122-A2148},
 volume = {36},
 number = {5},
 issn = {1064-8275},
 journal = {SIAM Journal on Scientific Computing},
 doi = {10.1137/130933381}
}

@article{Alexanderian.2018,
 author = {Alexanderian, Alen and Saibaba, Arvind K.},
 year = {2018},
 title = {Efficient D-Optimal Design of Experiments for Infinite-Dimensional Bayesian Linear Inverse Problems},
 pages = {A2956-A2985},
 volume = {40},
 number = {5},
 issn = {1064-8275},
 journal = {SIAM Journal on Scientific Computing},
 doi = {10.1137/17M115712X}
}

@article{Attia.2018,
 author = {Attia, Ahmed and Alexanderian, Alen and Saibaba, Arvind K.},
 year = {2018},
 title = {Goal-oriented optimal design of experiments for large-scale Bayesian linear inverse problems},
 pages = {095009},
 volume = {34},
 number = {9},
 journal = {Inverse Problems},
 doi = {10.1088/1361-6420/aad210}
}

@article{Biccari.2023,
 author = {Biccari, Umberto and Song, Yongcun and Yuan, Xiaoming and Zuazua, Enrique},
 year = {2023},
 title = {A two-stage numerical approach for the sparse initial source identification of a diffusion--advection equation},
 pages = {095003},
 volume = {39},
 number = {9},
 journal = {Inverse Problems},
 doi = {10.1088/1361-6420/ace548}
}

@inproceedings{Bonari.2024,
 author = {Bonari, Jacopo and K{\"u}hn, Lisa and von Danwitz, Max and Popp, Alexander},
 title = {Towards Real-Time Urban Physics Simulations with Digital Twins},
 pages = {18--25},
 publisher = {IEEE},
 isbn = {979-8-3315-2721-1},
 booktitle = {2024 28th International Symposium on Distributed Simulation and Real Time Applications (DS-RT)},
 year = {2024},
 doi = {10.1109/DS-RT62209.2024.00013}
}

@article{Boris.2002,
 author = {Boris, J.},
 year = {2002},
 title = {The threat of chemical and biological terrorism: preparing a response},
 pages = {22--32},
 volume = {4},
 number = {2},
 issn = {15219615},
 journal = {Computing in Science {\&} Engineering},
 doi = {10.1109/5992.988644}
}

@article {Boyer.2019,
    AUTHOR = {Boyer, Claire and Chambolle, Antonin and De Castro, Yohann and
              Duval, Vincent and de Gournay, Fr\'ed\'eric and Weiss, Pierre},
     TITLE = {On representer theorems and convex regularization},
   JOURNAL = {SIAM J. Optim.},
  FJOURNAL = {SIAM Journal on Optimization},
    VOLUME = {29},
      YEAR = {2019},
    NUMBER = {2},
     PAGES = {1260--1281},
      ISSN = {1052-6234,1095-7189},
   MRCLASS = {52A05 (46E27 49N45)},
  MRNUMBER = {3948246},
MRREVIEWER = {Qinghong\ Zhang},
       DOI = {10.1137/18M1200750},
}

@InProceedings{Scherzer.2009,
author="Scherzer, Otmar
and Walch, Birgit",
editor="Tai, Xue-Cheng
and M{\o}rken, Knut
and Lysaker, Marius
and Lie, Knut-Andreas",
title="Sparsity Regularization for Radon Measures",
booktitle="Scale Space and Variational Methods in Computer Vision",
year="2009",
publisher="Springer Berlin Heidelberg",
address="Berlin, Heidelberg",
pages="452--463",
doi={10.1007/978-3-642-02256-2_38}
}

@article{Bredies.2013,
 author = {Bredies, Kristian and Pikkarainen, Hanna Katriina},
 title = {Inverse problems in spaces of measures},
 fjournal = {European Series in Applied and Industrial Mathematics (ESAIM): Control, Optimization and Calculus of Variations},
 journal = {ESAIM, Control Optim. Calc. Var.},
 issn = {1292-8119},
 volume = {19},
 number = {1},
 pages = {190--218},
 year = {2013},
 language = {English},
 doi = {10.1051/cocv/2011205},
 zbMATH = {6144572},
 Zbl = {1266.65083}
}

@article{Chizat.2022,
 author = {Chizat, L{\'e}na{\"{\i}}c},
 title = {Sparse optimization on measures with over-parameterized gradient descent},
 fjournal = {Mathematical Programming. Series A. Series B},
 journal = {Math. Program.},
 issn = {0025-5610},
 volume = {194},
 number = {1-2 (A)},
 pages = {487--532},
 year = {2022},
 language = {English},
 doi = {10.1007/s10107-021-01636-z},
 zbMATH = {7550210},
 Zbl = {1494.90082}
}

@article{Denoyelle.2020,
 author = {Denoyelle, Quentin and Duval, Vincent and Peyr{\'e}, Gabriel and Soubies, Emmanuel},
 title = {The sliding {Frank}-{Wolfe} algorithm and its application to super-resolution microscopy},
 fjournal = {Inverse Problems},
 journal = {Inverse Probl.},
 issn = {0266-5611},
 volume = {36},
 number = {1},
 pages = {42},
 note = {Id/No 014001},
 year = {2020},
 language = {English},
 doi = {10.1088/1361-6420/ab2a29},
 zbMATH = {7153864},
 Zbl = {1434.65082}
}

@misc{Hnatiuk.2025,
 author = {Hnatiuk, Arsen and Walter, Daniel},
 title = {Lazifying point insertion algorithms in spaces of measures},
 year = {2025},
 doi={10.48550/arXiv.2508.03459}
}

@article {Duval.2015,
    AUTHOR = {Duval, Vincent and Peyr\'e, Gabriel},
     TITLE = {Exact support recovery for sparse spikes deconvolution},
   JOURNAL = {Found. Comput. Math.},
  FJOURNAL = {Foundations of Computational Mathematics. The Journal of the
              Society for the Foundations of Computational Mathematics},
    VOLUME = {15},
      YEAR = {2015},
    NUMBER = {5},
     PAGES = {1315--1355},
      ISSN = {1615-3375,1615-3383},
   MRCLASS = {65J20 (46E27 49K40 65R32)},
  MRNUMBER = {3394712},
MRREVIEWER = {Anton\ Suhadolc},
       DOI = {10.1007/s10208-014-9228-6},
}

@book{Boyd,
 author = {Boyd, Stephen P. and Vandenberghe, Lieven},
 year = {2004},
 title = {Convex optimization},
 address = {Cambridge and New York, NY and Port Melbourne},
 publisher = {{Cambridge University Press}},
 isbn = {978-0-521-83378-3},
 doi={10.1017/CBO9780511804441},
}

@article {Bredies.2025,
    AUTHOR = {Bredies, Kristian and Iglesias, Jos\'e{} A. and Walter,
              Daniel},
     TITLE = {On extremal points for some vectorial total variation
              seminorms},
   JOURNAL = {ESAIM Control Optim. Calc. Var.},
  FJOURNAL = {ESAIM. Control, Optimisation and Calculus of Variations},
    VOLUME = {31},
      YEAR = {2025},
     PAGES = {Paper No. 92, 35},
      ISSN = {1292-8119,1262-3377},
   MRCLASS = {46N10 (26B30 46A55 49Q20)},
  MRNUMBER = {4989255},
       DOI = {10.1051/cocv/2025078},
}

@article {Carioni.2020,
    AUTHOR = {Bredies, Kristian and Carioni, Marcello},
     TITLE = {Sparsity of solutions for variational inverse problems with
              finite-dimensional data},
   JOURNAL = {Calc. Var. Partial Differential Equations},
  FJOURNAL = {Calculus of Variations and Partial Differential Equations},
    VOLUME = {59},
      YEAR = {2020},
    NUMBER = {1},
     PAGES = {Paper No. 14, 26},
      ISSN = {0944-2669,1432-0835},
   MRCLASS = {49J27 (49N45 52A05)},
  MRNUMBER = {4040623},
MRREVIEWER = {Tamaz\ Tadumadze},
       DOI = {10.1007/s00526-019-1658-1},
}

@article{Brooks.1982,
 author = {Brooks, Alexander N. and Hughes, Thomas J.R.},
 year = {1982},
 title = {Streamline upwind/Petrov-Galerkin formulations for convection dominated flows with particular emphasis on the incompressible Navier-Stokes equations},
 pages = {199--259},
 volume = {32},
 number = {1-3},
 journal = {Computer Methods in Applied Mechanics and Engineering},
 doi = {10.1016/0045-7825(82)90071-8}
}

@article{Casas.2012,
 author = {Casas, Eduardo and Clason, Christian and Kunisch, Karl},
 year = {2012},
 title = {Approximation of Elliptic Control Problems in Measure Spaces with Sparse Solutions},
 pages = {1735--1752},
 volume = {50},
 number = {4},
 issn = {0363-0129},
 journal = {SIAM Journal on Control and Optimization},
 doi = {10.1137/110843216}
}

@article{Danwitz.2019,
 author = {von Danwitz, Max and Karyofylli, Violeta and Hosters, Norbert and Behr, Marek},
 year = {2019},
 title = {Simplex space--time meshes in compressible flow simulations},
 pages = {29--48},
 volume = {91},
 number = {1},
 issn = {0271-2091},
 journal = {International Journal for Numerical Methods in Fluids},
 doi = {10.1002/fld.4743}
}

@misc{Danwitz.2024,
 author = {von Danwitz, Max and Bonari, Jacopo and Franz, Philip and K{\"u}hn, Lisa and Mattuschka, Marco and Popp, Alexander},
 title = {Contaminant Dispersion Simulation in a Digital Twin Framework for Critical Infrastructure Protection},
 pages = {1--12},
 booktitle = {Towards Digital Twins for Infrastructures},
 doi = {10.23967/eccomas.2024.301}
}

@article{Danwitz.2023,
 author = {von Danwitz, Max and Voulis, Igor and Hosters, Norbert and Behr, Marek},
 year = {2023},
 title = {Time--continuous and time--discontinuous space--time finite elements for advection--diffusion problems},
 pages = {3117--3144},
 volume = {124},
 number = {14},
 journal = {International Journal for Numerical Methods in Engineering},
 doi = {10.1002/nme.7241}
}

@book{Elman.2014,
 author = {Elman, Howard C. and Silvester, David J. and Wathen, Andrew J.},
 year = {2014},
 title = {Finite elements and fast iterative solvers: with applications in incompressible fluid dynamics},
 address = {Oxford},
 edition = {Second Edition},
 publisher = {{Oxford University Press}},
 isbn = {9780199678808},
 series = {Numerical mathematics and scientific computation},
 doi={10.1093/acprof:oso/9780199678792.001.0001}
}

@article{Halko.2011,
 author = {Halko, N. and Martinsson, P. G. and Tropp, J. A.},
 year = {2011},
 title = {Finding Structure with Randomness: Probabilistic Algorithms for Constructing Approximate Matrix Decompositions},
 pages = {217--288},
 volume = {53},
 number = {2},
 issn = {0036-1445},
 journal = {SIAM Review},
 doi = {10.1137/090771806}
}

@article {Fernandez.2014,
    AUTHOR = {Cand\`es, Emmanuel J. and Fernandez-Granda, Carlos},
     TITLE = {Towards a mathematical theory of super-resolution},
   JOURNAL = {Comm. Pure Appl. Math.},
  FJOURNAL = {Communications on Pure and Applied Mathematics},
    VOLUME = {67},
      YEAR = {2014},
    NUMBER = {6},
     PAGES = {906--956},
      ISSN = {0010-3640,1097-0312},
   MRCLASS = {94A08 (65J22 65K05 94A12)},
  MRNUMBER = {3193963},
MRREVIEWER = {Joseph\ D.\ Lakey},
       DOI = {10.1002/cpa.21455},
}

@article{Huynh.2024,
 author = {Huynh, Phuoc-Truong and Pieper, Konstantin and Walter, Daniel},
 year = {2024},
 title = {Towards optimal sensor placement for inverse problems in spaces of measures},
 pages = {055007},
 volume = {40},
 number = {5},
 journal = {Inverse Problems},
 doi = {10.1088/1361-6420/ad2cf8}
}

@article{Key.2023,
 author = {Key, Fabian and von Danwitz, Max and Ballarin, Francesco and Rozza, Gianluigi},
 year = {2023},
 title = {Model order reduction for deforming domain problems in a time--continuous space--time setting},
 pages = {5125--5150},
 volume = {124},
 number = {23},
 journal = {International Journal for Numerical Methods in Engineering},
 doi = {10.1002/nme.7342}
}

@article{neitzel19,
 author = {Neitzel, Ira and Pieper, Konstantin and Vexler, Boris and Walter, Daniel},
 title = {A sparse control approach to optimal sensor placement in {PDE}-constrained parameter estimation problems},
 fjournal = {Numerische Mathematik},
 journal = {Numer. Math.},
 issn = {0029-599X},
 volume = {143},
 number = {4},
 pages = {943--984},
 year = {2019},
 language = {English},
 doi = {10.1007/s00211-019-01073-3},
 zbMATH = {7127879},
 Zbl = {1423.35452}
}

@article{Milzarek.2014,
 author = {Milzarek, Andre and Ulbrich, Michael},
 year = {2014},
 title = {A Semismooth Newton Method with Multidimensional Filter Globalization for $l_1$-Optimization},
 pages = {298--333},
 volume = {24},
 number = {1},
 issn = {1052-6234},
 journal = {SIAM Journal on Optimization},
 doi = {10.1137/120892167}
}

@article{Monge.2020,
 author = {Monge, Azahar and Zuazua, Enrique},
 year = {2020},
 title = {Sparse source identification of linear diffusion--advection equations by adjoint methods},
 pages = {104801},
 volume = {145},
 journal = {Systems {\&} Control Letters},
 doi = {10.1016/j.sysconle.2020.104801}
}

@incollection{Patnaik.2012,
 author = {Patnaik, Gopal and Boris, Jay P. and Grinstein, Fernando F. and Iselin, John P. and Hertwig, Denise},
 title = {Large Scale Urban Simulations with FCT},
 pages = {91--117},
 publisher = {{Springer Netherlands}},
 isbn = {978-94-007-4037-2 },
 series = {Scientific Computation},
 editor = {Kuzmin, Dmitri and L{\"o}hner, Rainald and Turek, Stefan},
 booktitle = {Flux-Corrected Transport},
 year = {2012},
 address = {Dordrecht},
 doi = {10.1007/978-94-007-4038-9_4 }
}

@misc{Petra.2011,
 author = {Petra, Noei and Stadler, Georg},
 date = {2011},
 title = {Model Variational Inverse Problems Governed by Partial Differential Equations},
 howpublished = {\url{https://math.nyu.edu/~stadler/papers/PetraStadler11.pdf}},
 address = {The University of Texas},
 institution = {{The Institute for Computational Engineering and Sciences}}
}

@article{Pieper.2021,
 author = {Pieper, Konstantin and Walter, Daniel},
 year = {2021},
 title = {Linear convergence of accelerated conditional gradient algorithms in spaces of measures},
 pages = {38},
 volume = {27},
 issn = {1292-8119},
 journal = {ESAIM: Control, Optimisation and Calculus of Variations},
 doi = {10.1051/cocv/2021042 }
}

@article{Steihaug.1983,
 author = {Steihaug, Trond},
 year = {1983},
 title = {Local and superlinear convergence for truncated iterated projections methods},
 pages = {176--190},
 volume = {27},
 number = {2},
 issn = {0025-5610},
 journal = {Mathematical Programming},
 doi = {10.1007/BF02591944}
}

@article{Villa.2021,
 author = {Villa, Umberto and Petra, Noemi and Ghattas, Omar},
 year = {2021},
 title = {hIPPYlib},
 pages = {1--34},
 volume = {47},
 number = {2},
 journal = {ACM Transactions on Mathematical Software},
 doi = {10.1145/3428447}
}

@incollection{Wogrin.2023,
 author = {Wogrin, Sonja and Singh, Arjun and Allaire, Douglas and Ghattas, Omar and Willcox, Karen},
 title = {From Data to Decisions: A Real-Time Measurement--Inversion--Prediction--Steering Framework for Hazardous Events and Health Monitoring},
 pages = {195--227},
 publisher = {{Springer International Publishing}},
 isbn = {978-3-031-27985-0},
 editor = {Darema, Frederica and Blasch, Erik P. and Ravela, Sai and Aved, Alex J.},
 booktitle = {Handbook of Dynamic Data Driven Applications Systems},
 year = {2023},
 address = {Cham},
 doi = {10.1007/978-3-031-27986-7_8}
}

@article{Wu.2023,
 author = {Wu, Keyi and Chen, Peng and Ghattas, Omar},
 year = {2023},
 title = {An Offline-Online Decomposition Method for Efficient Linear Bayesian Goal-Oriented Optimal Experimental Design: Application to Optimal Sensor Placement},
 pages = {B57-B77},
 volume = {45},
 number = {1},
 issn = {1064-8275},
 journal = {SIAM Journal on Scientific Computing},
 doi = {10.1137/21M1466542}
}

@article{Daon.2018,
 author = {Daon, Yair and Stadler, Georg},
 year = {2018},
 title = {Mitigating the influence of the boundary on PDE-based covariance operators},
 pages = {1083--1102},
 volume = {12},
 number = {5},
 issn = {1930-8345},
 journal = {Inverse Problems {\&} Imaging},
 doi = {10.3934/ipi.2018045}
}

@article{Casas.2019,
 author = {Casas, E. and Kunisch, K.},
 year = {2019},
 title = {Using sparse control methods to identify sources in linear diffusion-convection equations},
 pages = {114002},
 volume = {35},
 number = {11},
 issn = {0266-5611},
 journal = {Inverse Problems},
 doi = {10.1088/1361-6420/ab331c}
}

@article{Leykekhman.2020,
 author = {Leykekhman, Dmitriy and Vexler, Boris and Walter, Daniel},
 year = {2020},
 title = {Numerical analysis of sparse initial data identification for parabolic problems},
 pages = {1139--1180},
 volume = {54},
 number = {4},
 issn = {0764-583X},
 journal = {ESAIM: Mathematical Modelling and Numerical Analysis},
 doi = {10.1051/m2an/2019083}
}

@article{Gorelick.1983,
 author = {Gorelick, Steven M. and Evans, Barbara and Remson, Irwin},
 year = {1983},
 title = {Identifying sources of groundwater pollution: An optimization approach},
 pages = {779--790},
 volume = {19},
 number = {3},
 journal = {Water Resources Research},
 doi = {10.1029/WR019i003p00779}
}

@article{Tsai.2014,
 author = {Tsai, Richard and Osher, Stanley and Li, Yingying},
 year = {2014},
 title = {Heat source identification based on $l_1$ constrained minimization},
 pages = {199--221},
 volume = {8},
 number = {1},
 journal = {Inverse Problems and Imaging},
 doi = {10.3934/ipi.2014.8.199}
}

@article{Elman.2020,
 author = {Elman, Howard C. and Su, Tengfei},
 year = {2020},
 title = {A low-rank solver for the stochastic unsteady Navier--Stokes problem},
 pages = {112948},
 volume = {364},
 issn = {0045-7825},
 journal = {Computer Methods in Applied Mechanics and Engineering},
 doi = {10.1016/j.cma.2020.112948}
}

@inproceedings{Ruiz.2024,
 author = {Ruiz, Victor Prieto and Hinsen, Patrick and Wiedemann, Thomas and Shutin, Dmitriy and Christof, Constantin},
 title = {Gas Source Localization Using Physics-Guided Neural Networks},
 pages = {1--3},
 publisher = {IEEE},
 isbn = {979-8-3503-4865-1},
 booktitle = {2024 IEEE International Symposium on Olfaction and Electronic Nose (ISOEN)},
 year = {2024},
 doi = {10.1109/ISOEN61239.2024.10556061}
}

@inproceedings{Hinsen.2024,
 author = {Hinsen, Patrick and Wiedemann, Thomas and Ruiz, Victor Scott Prieto and Shutin, Dmitriy and Lilienthal, Achim J.},
 title = {Experimental Study of Gas Propagation: Parameter Identification and Analysis in a Wind Tunnel},
 pages = {1--3},
 publisher = {IEEE},
 isbn = {979-8-3503-4865-1},
 booktitle = {2024 IEEE International Symposium on Olfaction and Electronic Nose (ISOEN)},
 year = {2024},
 doi = {10.1109/ISOEN61239.2024.10556306}
}

@inproceedings{Wiedemann.2024,
 author = {Wiedemann, Thomas and Hinsen, Patrick and Ruiz, Victor Prieto and Shutin, Dmitriy and Lilienthal, Achim J.},
 title = {Domain Knowledge Assisted Gas Tomography},
 pages = {1--3},
 publisher = {IEEE},
 isbn = {979-8-3503-4865-1},
 booktitle = {2024 IEEE International Symposium on Olfaction and Electronic Nose (ISOEN)},
 year = {2024},
 doi = {10.1109/ISOEN61239.2024.10556226}
}

@article{Arystanbekova.2004,
 author = {Arystanbekova, N.Kh.},
 year = {2004},
 title = {Application of Gaussian plume models for air pollution simulation at instantaneous emissions},
 pages = {451--458},
 volume = {67},
 number = {4-5},
 journal = {Mathematics and Computers in Simulation},
 doi = {10.1016/j.matcom.2004.06.023}
}

@article{Ulfah.2018,
 author = {Ulfah, S. and Awalludin, S. A. and Wahidin},
 year = {2018},
 title = {Advection-diffusion model for the simulation of air pollution distribution from a point source emission},
 pages = {012067},
 volume = {948},
 journal = {Journal of Physics: Conference Series},
 doi = {10.1088/1742-6596/948/1/012067}
}

@article{Albani.2021,
 author = {Albani, Roseane A. S. and Albani, Vinicius V. L. and Migon, H{\'e}lio S. and {Silva Neto}, Ant{\^o}nio J.},
 year = {2021},
 title = {Uncertainty quantification and atmospheric source estimation with a discrepancy-based and a state-dependent adaptative MCMC},
 pages = {118039},
 volume = {290},
 journal = {Environmental Pollution},
 doi = {10.1016/j.envpol.2021.118039}
}

@article{Hutchinson.2017,
 author = {Hutchinson, Michael and Oh, Hyondong and Chen, Wen-Hua},
 year = {2017},
 title = {A review of source term estimation methods for atmospheric dispersion events using static or mobile sensors},
 pages = {130--148},
 volume = {36},
 journal = {Information Fusion},
 doi = {10.1016/j.inffus.2016.11.010}
}

@article{Holmes.2006,
 author = {Holmes, N. S. and Morawska, L.},
 year = {2006},
 title = {A review of dispersion modelling and its application to the dispersion of particles: An overview of different dispersion models available},
 pages = {5902--5928},
 volume = {40},
 number = {30},
 journal = {Atmospheric Environment},
 doi = {10.1016/j.atmosenv.2006.06.003}
}

@article{Dinkel.2024,
 author = {Dinkel, Maximilian and Geitner, Carolin M. and {Robalo Rei}, Gil and Nitzler, Jonas and Wall, Wolfgang A.},
 year = {2024},
 title = {Solving Bayesian inverse problems with expensive likelihoods using constrained Gaussian processes and active learning},
 pages = {095008},
 volume = {40},
 number = {9},
 issn = {0266-5611},
 journal = {Inverse Problems},
 doi = {10.1088/1361-6420/ad5eb4}
}

@article{Nitzler.2022,
 author = {Nitzler, Jonas and Biehler, Jonas and Fehn, Niklas and Koutsourelakis, Phaedon-Stelios and Wall, Wolfgang A.},
 year = {2022},
 title = {A generalized probabilistic learning approach for multi-fidelity uncertainty quantification in complex physical simulations},
 pages = {115600},
 volume = {400},
 issn = {0045-7825},
 journal = {Computer Methods in Applied Mechanics and Engineering},
 doi = {10.1016/j.cma.2022.115600}
}

@article{Maronga.2020,
 author = {Maronga, Bj{\"o}rn and Banzhaf, Sabine and Burmeister, Cornelia and Esch, Thomas and Forkel, Renate and Fr{\"o}hlich, Dominik and Fuka, Vladimir and Gehrke, Katrin Frieda and Geleti{\v{c}}, Jan and Giersch, Sebastian and Gronemeier, Tobias and Gro{\ss}, G{\"u}nter and Heldens, Wieke et al.},
 year = {2020},
 title = {Overview of the PALM model system 6.0},
 pages = {1335--1372},
 volume = {13},
 number = {3},
 journal = {Geoscientific Model Development},
 doi = {10.5194/gmd-13-1335-2020}
}

@article{Chen.2025,
 author = {Chen, Yuanqing and Wang, Ding and Feng, Dachuan and Tian, Geng and Gupta, Vikrant and Cao, Renjing and Wan, Minping and Chen, Shiyi},
 year = {2025},
 title = {Three-dimensional spatiotemporal wind field reconstruction based on LiDAR and multi-scale PINN},
 pages = {124577},
 volume = {377},
 journal = {Applied Energy},
 doi = {10.1016/j.apenergy.2024.124577}
}

@article{Towers.2016,
 author = {Towers, P. and Jones, B. Ll.},
 year = {2016},
 title = {Real-time wind field reconstruction from LiDAR measurements using a dynamic wind model and state estimation},
 pages = {133--150},
 volume = {19},
 number = {1},
 journal = {Wind Energy},
 doi = {10.1002/we.1824}
}

@article{Sun.2022,
 author = {Sun, Shilin and Wang, Tianyang and Chu, Fulei and Tan, Jianxin},
 year = {2022},
 title = {Acoustic source identification using an off-grid and sparsity-based method for sound field reconstruction},
 pages = {108869},
 volume = {170},
 issn = {08883270},
 journal = {Mechanical Systems and Signal Processing},
 doi = {10.1016/j.ymssp.2022.108869}
}

@article{Jerez.2021,
 author = {Jerez, D. J. and Jensen, H. A. and Beer, M. and Broggi, M.},
 year = {2021},
 title = {Contaminant source identification in water distribution networks: A Bayesian framework},
 pages = {107834},
 volume = {159},
 issn = {08883270},
 journal = {Mechanical Systems and Signal Processing},
 doi = {10.1016/j.ymssp.2021.107834}
}

@article{Scholl.2023,
 author = {Scholl, Manuel and Passek, Matthias and Lainer, Mirjam and Taddei, Francesca and Schneider, Felix and M{\"u}ller, Gerhard},
 year = {2023},
 title = {Data-based localisation methods using simulated data with application to small-scale structures},
 pages = {2411--2433},
 volume = {93},
 number = {6},
 journal = {Archive of Applied Mechanics},
 doi = {10.1007/s00419-023-02388-2}
}

@article{Khodayimehr.2019,
 author = {Khodayi-mehr, Reza and Aquino, Wilkins and Zavlanos, Michael M.},
 year = {2019},
 title = {Model-Based Active Source Identification in Complex Environments},
 pages = {633--652},
 volume = {35},
 number = {3},
 journal = {IEEE Transactions on Robotics},
 doi = {10.1109/TRO.2019.2894039}
}

@article{Mazzieri.2022,
 author = {Mazzieri, Ilario and Muhr, Markus and Stupazzini, Marco and Wohlmuth, Barbara},
 year = {2022},
 title = {Elasto-acoustic modeling and simulation for the seismic response of structures: The case of the Tahtal{\i} dam in the 2020 Izmir earthquake},
 pages = {111411},
 volume = {466},
 issn = {00219991},
 journal = {Journal of Computational Physics},
 doi = {10.1016/j.jcp.2022.111411}
}

@article{Gjerde.2022,
 author = {Gjerde, Ingeborg and Scott, L.},
 year = {2022},
 title = {Nitsche's method for Navier--Stokes equations with slip boundary conditions},
 pages = {597--622},
 volume = {91},
 number = {334},
 journal = {Mathematics of Computation},
 doi = {10.1090/mcom/3682}
}

@article{Schneider.2025,
 author = {Schneider, Moritz and Halekotte, Lukas and Comes, Tina and Lichte, Daniel and Fiedrich, Frank},
 year = {2025},
 title = {Emergency Response Inference Mapping (ERIMap): A Bayesian network-based method for dynamic observation processing},
 pages = {110640},
 volume = {255},
 journal = {Reliability Engineering {\&} System Safety},
 doi = {10.1016/j.ress.2024.110640}
}

@article{Brusca.2016,
 author = {Brusca, S. and Famoso, F. and Lanzafame, R. and Mauro, S. and Garrano, A. Marino Cugno and Monforte, P.},
 year = {2016},
 title = {Theoretical and Experimental Study of Gaussian Plume Model in Small Scale System},
 pages = {58--65},
 volume = {101},
 journal = {Energy Procedia},
 doi = {10.1016/j.egypro.2016.11.008}
}

@phdthesis{walter.dissertation,
	author = {Walter, Daniel},
	title = {On sparse sensor placement for parameter identification problems with partial differential equations},
	year = {2019},
	school = {Technische Universität München},
	pages = {337},
	language = {en},
}

@article{Stuart.2010,
 author = {Stuart, A. M.},
 year = {2010},
 title = {Inverse problems: A Bayesian perspective},
 pages = {451--559},
 volume = {19},
 journal = {Acta Numerica},
 doi = {10.1017/S0962492910000061}
}

@misc{Mattuschka.03.07.2025,
 author = {Mattuschka, Marco and Lan, Noah An der and von Danwitz, Max and Wolff, Daniel and Popp, Alexander},
 date = {03.07.2025},
 title = {Goal-oriented optimal sensor placement for PDE-constrained inverse  problems in crisis management},
 doi={10.48550/arXiv.2507.02500}
}

@book{Friedman.2008,
 author = {Friedman, Avner},
 year = {2008},
 title = {Partial differential equations of parabolic type},
 address = {Mineola, New York},
 edition = {Dover ed.},
 publisher = {{Dover Publications}},
 isbn = {0486466256}
}

@article{Givoli.2014,
 author = {Givoli, Dan},
 year = {2014},
 title = {Time Reversal as a Computational Tool in Acoustics and Elastodynamics},
 pages = {1430001},
 volume = {22},
 number = {03},
 journal = {Journal of Computational Acoustics},
 doi = {10.1142/S0218396X14300011}
}



\end{document}